\documentclass[12pt]{spieman}  
\RequirePackage{fix-cm}
\usepackage{cancel}
\usepackage{amsmath,amsfonts,amssymb}
\usepackage{graphicx}
\usepackage[colorlinks=true,allcolors=blue]{hyperref}
\usepackage{tocloft}

\usepackage[monochrome]{xcolor}              
\usepackage[nomarkers]{endfloat} 
\usepackage{wrapfig}             

\usepackage{algorithm}
\usepackage{algpseudocode}
\algnewcommand{\LineComment}[1]{\State  \(\triangleright\) #1 \hfill~}

\title{ Three-dimensional multiscale discrete Radon and John transforms}

\author[a,*]{Jos\'e G. Marichal-Hern\'andez}
\author[a,b]{\'Oscar G\'omez-C\'ardenes}
\author[a]{Fernando Rosa}
\author[c]{Do Hyung Kim}
\author[a,b]{Jos\'e M. Rodr\'iguez-Ramos}
\affil[a]{Universidad de La Laguna, Industrial Engineering Department, ESIT, La Laguna, Spain, 38200.}
\affil[b]{Wooptix S.L., Av. Trinidad, 61, La Laguna, Spain, 38204.}
\affil[c]{Electronics and Telecommunications Research Institute, CG/Vision Research Group, 128 Gajeong-ro, Yuseong-gu, Daejeon, Korea.}

\begin{document} 
\maketitle
Note: This document is a preprint. The published manuscript can be found at \url{https://doi.org/10.1117/1.OE.59.9.093104}. © Copyright 2022 Society of Photo-Optical Instrumentation Engineers (SPIE). One print or electronic copy may be made for personal use only. Systematic reproduction and distribution, duplication of any material in this paper for a fee or for commercial purposes, or modification of the content of the paper are prohibited. 

\begin{abstract}
Two algorithms are introduced for the computation of discrete integral transforms with a multiscale approach operating in discrete three-dimensional (3D) volumes while considering its real-time implementation.

The first algorithm, referred to as 3D discrete Radon transform (DRT) of planes, will compute the summation set of values lying in discrete planes in a cube that imitates, in discrete data, the integrals on two-dimensional planes in a 3D volume similar to the continuous Radon transform. The normals of these planes, equispaced in ascents, cover a quadrilateralized hemisphere and comprise 12 dodecants.

The second proposed algorithm, referred to as the 3D discrete John transform (DJT) of lines, will sum elements lying on discrete 3D lines while imitating the behavior of the John or X-ray continuous transform on 3D volumes.

These discrete integral transforms do not perform interpolation on input or intermediate data, and they can be computed using only integer arithmetics with linearithmic complexity; thus, outperforming the methods based on the Fourier slice-projection theorem for real-time applications. 

We briefly prove that these transforms have fast inversion algorithms that are exact for discrete inputs.
\end{abstract}

{\noindent \footnotesize\textbf{*}Jos\'e G. Marichal-Hern\'andez,  \linkable{jmariher@ull.edu.es} }


\section{INTRODUCTION}
\label{sec:intro}  
\subsection{Radon and John transforms on 3D}
The Radon transform \cite{Radon17} of a function on $\mathfrak{R}^n$ evaluates the integral of the function over hyperplanes, which are subspaces of dimension $n-1$, that is, line integrals in a two-dimensional (2D) space or plane integrals for a three-dimensional (3D) space.
The John (or X-ray) transform \cite{john1938}, calculates line integrals independently in which an $n-$dimensional space is being considered, and it coincides with the Radon transform for $n=2$.

These integral transforms are of significant interest in tomography \cite{kak2001principles}, and when lines or planes should be evaluated or detected \cite{surveyHough, Geng:18}. We are interested in their numerical computation on discrete data, and accordingly, these transforms can be promptly computed using the Fourier slice-projection theorem and thereafter the trigonometric polynomials and fractional Fourier transforms. This method has been established for the Radon and John transforms\cite{3DradonFourier, 3DXrayFourier}, and the use of the Fourier slice-projection theorem to solve 3D Radon transforms can be traced back to the origins of computer-assisted tomography \cite{Sheep80}. 

We propose to solve similar problems using the founding idea behind multiscale discrete Radon transform (DRT) as a starting point, proposed by several authors \cite{GotzDruckmuller, Brady, Brandt}, for 2D in the late 1990s. The multiscale approach relies on the recursive approximation of a line by the union of two line segments; and consequently, the 2D Radon transform on a dyadic square can be computed by combining recursively the Radon transforms of its four dyadic subsquares.

Similarly, the 3D John and Radon transforms on a dyadic cube can be computed by combining the transforms of its eight dyadic subcubes.

Donoho \emph{et al.}\cite{Donoho02fastx-ray} faced the same problem regarding computing beamlets when employing dyadic decomposition of cubes. This study preconized the existence of the algorithm for the 3D discrete John transform (DJT) that we are fully describing now, in Section \emph{4.2 Two-scale recursion}. However, the algorithm was discarded in favor of a shearing method based on fractional Fourier transforms, which is similar to what Averbuch and Shkolnisky proposed. Both groups of researchers later jointly refined the Radon theory in 2D \cite{FSS, AFFI, AFFII}, based on the fractional Fourier transform.
More recently, an exact and non-iterative inversion of the 3D DRT of planes based on the Fourier transform has been proposed \cite{inverse3DFourierRadon}, by mixing the pseudo-polar Fourier transform and multilevel Toeplitz operators\cite{Rauth}.

There are not several studies that solve a 3D DRT without Fourier transforms. Among those, Levi and Efros\cite{SHAS} have proposed an alternative method for 2D and 3D Radon transforms that interpolate data exploiting the periodicity of weights of line-grid cuts. The improved interpolation of the method is gained by being slower than what we will propose here.


The crude approach to interpolation of multiscale methods has been disregarded in the past \cite{FSS} owing to its geometrical inconsistencies and properties of the continuum transform that are not inherited in the discrete transform. However, we still believe this method also poses several computational advantages: faster computation, parallel prone nature, elusion of trigonometrics and real numbers, and an equally fast inversion algorithm, as demonstrated by Press for two dimensions\cite{Press}.


\subsection{Multiscale discrete Radon transform}
The DRT was originally designed to compute all the integrals across picture elements located on a discrete line that touches at least one element on a picture of size $N \times N$ while projecting on a semi circumference around it. An algorithm exists that solves this problem with linearithmic complexity, $O(N^2 \log N)$, for a quadrant from 0 to $45^\circ$. Additionally, by combining four runs of the algorithm on four mirrored versions of the input, a sort of M{\"o}bius band is obtained comprising the entire set of line integrals that occur when projecting $180^\circ$ around the image; a discrete version of the \emph{sinogram} in the continuum devised by Radon\cite{Radon17}.

\begin{figure}[h]
\centering
\includegraphics[width=.6\textwidth]{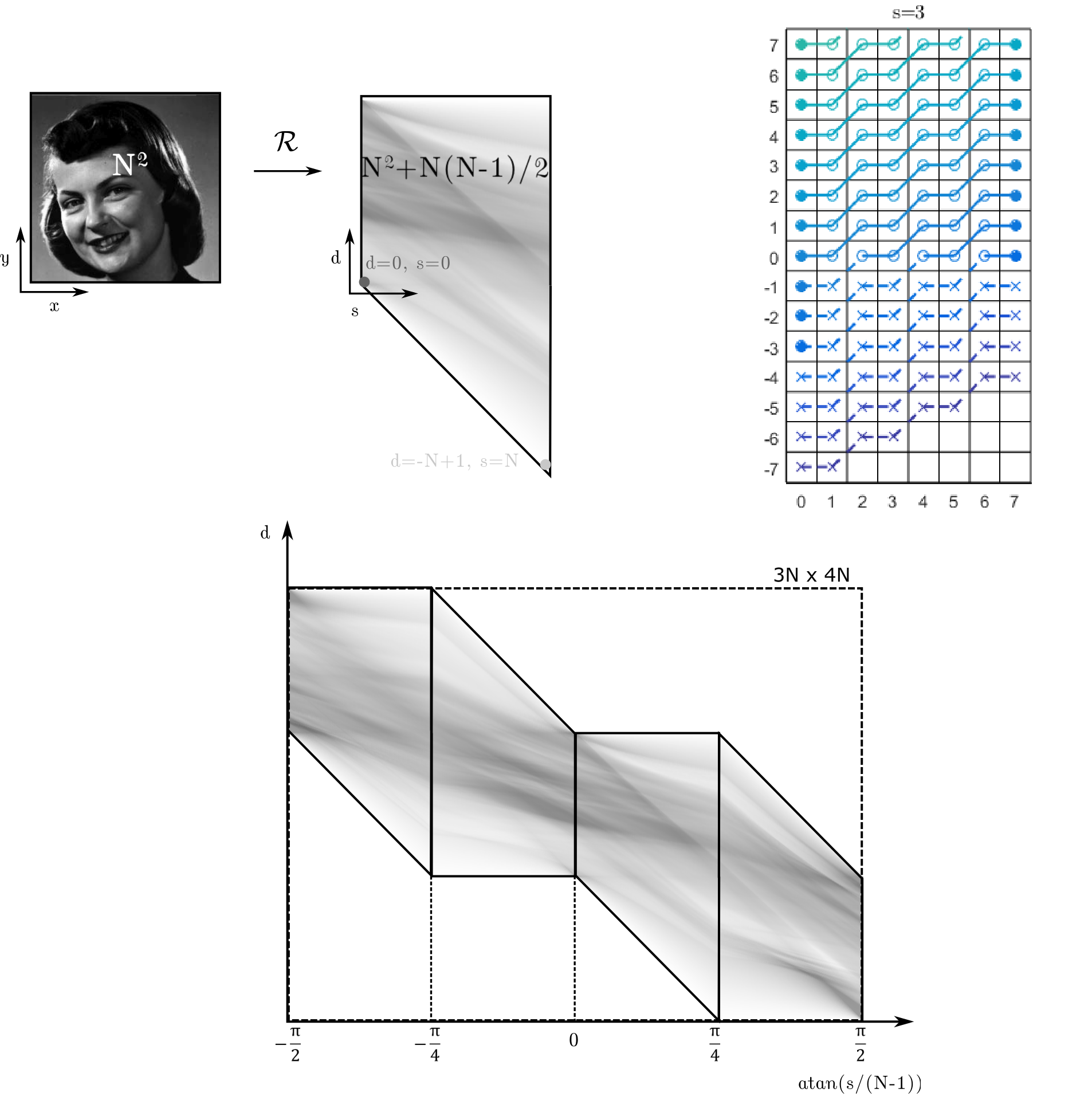}
\caption{Top left: Depiction of the relative sizes and shapes of an input and its 2D DRT for the quadrant covering $0^\circ$ to $45^\circ$. Top right: Set of displacements through an $8\times8$ domain for digital lines with slope 3.
Bottom: DRT of the four quadrants covering $180^\circ$ merged together.}
\label{fig:ConventionalShears}
\end{figure}

Figure \ref{fig:ConventionalShears} shows a grayscale image of size $N\times N$ and its DRT where pixels lying in the same discrete lines have been summed together. The lines are of the form $y = \frac{s}{N-1} \cdot x - d$, where $s \in \{0..N-1\}$ denotes the \emph{s}lope or ascent and $d \in \{-N+1..N-1\}$ denotes the \emph{d}isplacement or intercept. Moreover, a depiction of discrete lines is shown in an $8 \times 8$ domain for the slope $s=3$. This variety of names occurs due to the adoption of the algebraic slope-intercept formulation of lines and the adherence to the previous authors' notation, simultaneously. Henceforth, we will denote the parameters on a Radon domain as $s$ and $d$ for slope and displacement, respectively.

\subsection{Main contributions}
\begin{figure}[h]
\centering
\includegraphics[width=.5\textwidth]{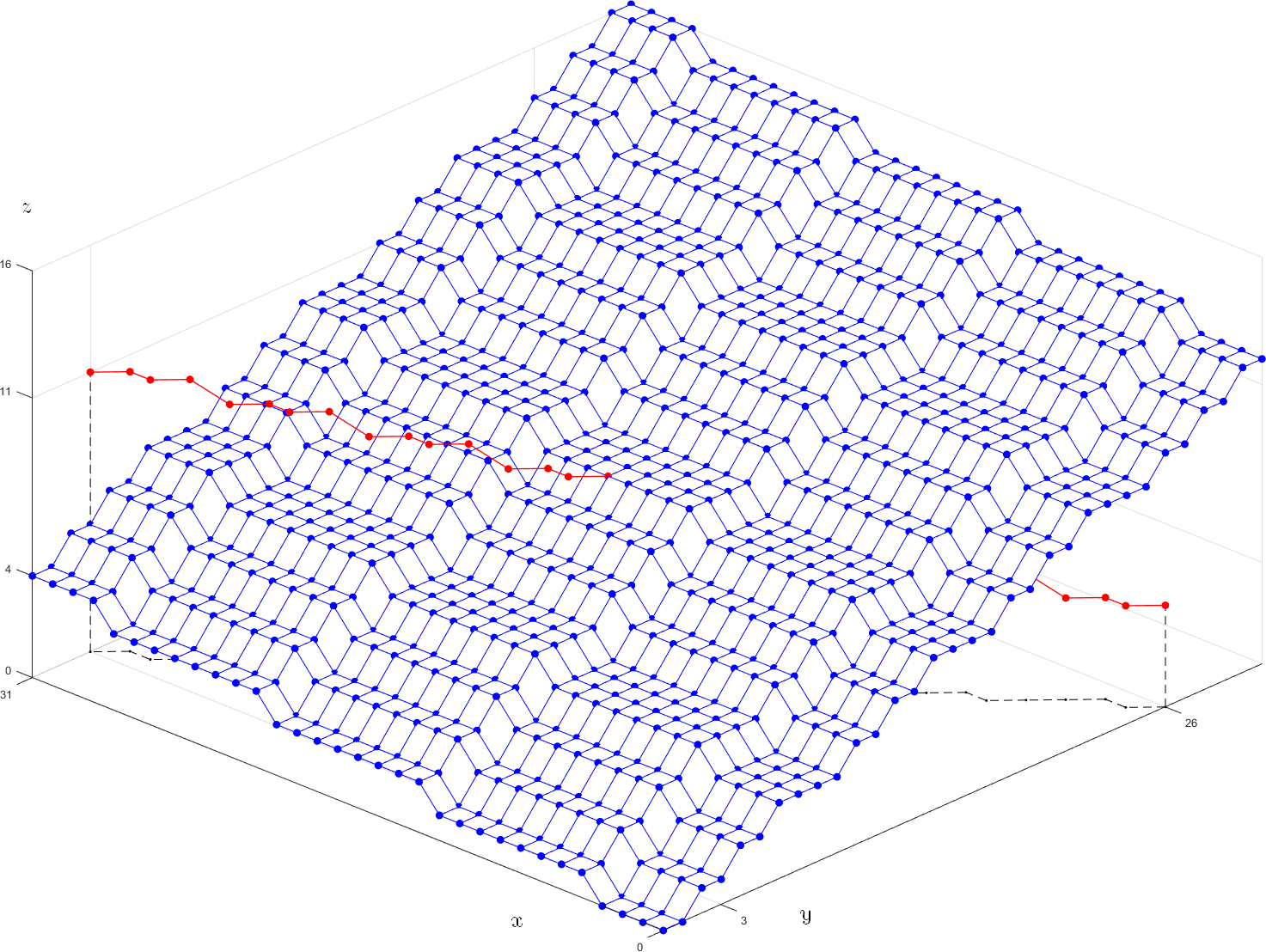} 
\caption{Discrete plane $z = l^5_{4}(\mathbf{x}) + l^5_{12}(\mathbf{y})$, and ray $(x,y,z) = ( \lambda(\mathbf{u}), \, 26 - l^5_{23}(\mathbf{u}), \, 4 + l^5_{7}(\mathbf{u})\, )$ }
\label{fig:planeLine3D}
\end{figure}
The main contributions of this study are as follows:
\begin{enumerate}
    \item  The 2D loose discrete line definition can be used to construct 3D discrete lines and planes, as shown in figure \ref{fig:planeLine3D}.
    \item  The discrete lines and planes are such that by recursion an algorithm can be created to compute the sets, which are complete in terms of slopes and intercepts, of summations of discrete data across them with linearithmic complexity.
    \item The forward algorithms that occur for planes and rays have equally fast backprojection algorithms.
    \item The pairs of forward and backprojection transforms can be combined to create a fast, although iterative, backward transform.
\end{enumerate} 

\subsection{Founding idea on 3D multiscale transforms}


\begin{figure}[h]
\centering
\includegraphics[width=.5\textwidth]{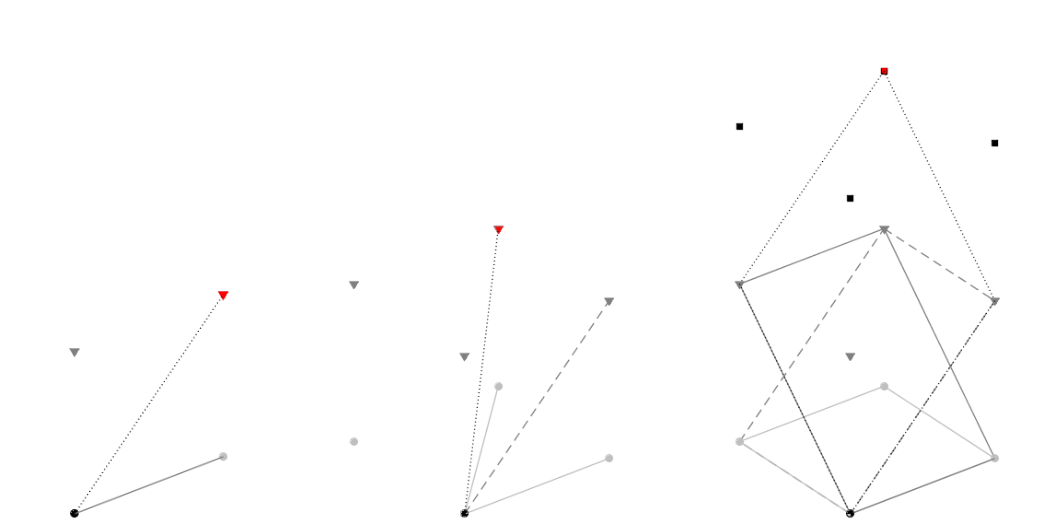}
\caption{Left: Two unique alternatives for a discrete line segment crossing the black dot in 2D at the lowest scale. Middle: Alternatives for a 3D line crossing through the black dot at the lowest scale.  Right: Four unique alternatives for a discrete plane patch crossing the black dot in 3D at the lowest scale.}
\label{fig:schemes}
\end{figure}

Figure \ref{fig:schemes}, adapted from Marichal \emph{et al.}\cite{Marichal}, depicts the foundation of discrete lines at the lowest scale in 2D, on the left, in which every discrete line crossing the black circle, with coordinates $\{2i, j\}$, are shown. The foundation will comprise each of the two line segments, the one joining the aforementioned point with the point at $\{2i+1, j\}$ or $\{2i+1, j+1\}$. The first stage of the partial transform of conventional 2D DRT involves adding the two combinations, $\forall i \in \{0..\frac{N}{2}-1\}, \forall j \in \{0..N-1\}$, and saving them as pertaining to slope 0 or 1 at the lowest scale in that subdomain of the problem. With these data, the algorithm can now operate at a greater scale by adding sums of two-point line segments from the previous stage to create every sum on a four-point line segment combination that produces slopes 0, $\frac{1}{3}$, $\frac{2}{3}$, and $\frac{3}{3}$. The method can operate again to create every possible sum of eight-point line segments, with slopes from 0 to $\frac{7}{7}$, etc. Then, after $\log_2(N)$ of these stages, the 2D DRT computation of a quadrant with an image sized $N \times N$ is completed.

This method applies in 3D volumes when dealing with sums of \emph{voxels} pertaining to every possible discrete line or plane that crosses through a cube.

In the case of planes, to the right of figure \ref{fig:schemes}, two slopes characterize the normal vector to the planes being considered, one slope running across the $x$-axis, referred to as $s_1$, and the other slope on the $y$-axis, $s_2$, in the basic \emph{quadrant} configuration. Therefore, four combinations emerge at the lowest scale.

For 3D lines, depicted in the middle of figure \ref{fig:schemes}, there are again two slopes while the algorithm advances from a generic point at $\{2i, j, k\}$ to the next position at $\{2i+1, j+s_1, k+s_2\}$ at the basic level.


To the best of our knowledge, the following have hitherto not been studied: the mapping equations from a lower stage to a higher stage for each of these 3D transforms; the number of the planes or lines that comprise the entire set of hyperplanes on the cube; and the manner in which the basic version of the algorithm is to be be applied to mirrored versions of the input to achieve this completeness.

\subsection{Structure of this paper}
We start by discussing the multiscale approach to DRT, followed by explaining the details of 3D DRT of planes in section \ref{sec:planes3D}. 
Then, we define the 3D DJT of lines in section \ref{sec:lines3D}.
Subsequently, we briefly comment on the invertibility of both algorithms in section \ref{sec:inversion}. The details and execution times of implementation are provided in section \ref{sec:implementationDetails}.
We will outline some future scope for study in section \ref{sec:conclusions} to conclude the paper.

\section{CONVENTIONAL DISCRETE 2D RADON TRANSFORM}
\label{sec:radon}
The 2D Radon transform allows a 2D signal, $f(\cdot)$, to be described in the continuum in terms of the integrals along lines parameterized using an angle and a displacement, $(\theta, \rho)$, instead of single values accessed by their horizontal and vertical pair of Cartesian coordinates $(x, y)$.
That is, \begin{equation} \mathfrak{R}f(\theta, \rho) = \iint f(x,y) \;
\delta(x\;cos\;\theta+y\;sin\;\theta-\rho) \; \mathrm{d} x \; \mathrm{d} y, \end{equation}
or equivalently using the absolute \emph{s}lope and \emph{d}isplacement form, with $|s| < 1$,
\begin{equation}
\begin{split}
\mathfrak{R}_{|\theta| \leq \frac{\pi}{2}}f(s, d) = \int f(u,\; u \; s + d) \;  \mathrm{d} u, \\
\mathfrak{R}_{|\theta| \geq \frac{\pi}{2}}f(s, d) = \int f(u \; s + d,\; u) \;  \mathrm{d} u.
\label{eq:radon2D} 
\end{split}
\end{equation}

The computer calculation of the Radon transform has to account for the discretization of data. If we have regularly spaced discrete data, the data will be available only for a finite number of samples normally accessed with integer indexes. The problem occurs when the continuous definition of the line integral makes it necessary to evaluate the function at positions where there is no sample available and interpolation is required.

By using the pseudo-polar Fourier transform \cite{AFFI} --a variation of the fast Fourier transform (FFT) that operates on a grid of concentric squares--, a DRT that is algebraically exact, invertible, fast, \cite{AFFII} and able to be generalized to 3D \cite{FSS,3DradonFourier} can be designed. However, it is based on Fourier transforms, which is best avoided because even if all those families of transforms are \emph{asymptotically} equal, $O(N^2 \log N)$ --in terms of computational complexity--, their runtime complexity includes multipliers due to resampling that can cause a significant difference when dealing with big data. 
Moreover, if the problem that should be solved has discrete inputs, the multiscale Radon transform is more accurate than other discretizations of Radon transforms \cite{kingston}.

\subsection{Forward discrete Radon transform}
Almost simultaneously, several authors\cite{GotzDruckmuller, Brady, Brandt} proposed a divide-and-conquer approach reminiscent of FFT in that it solves the problem at reduced scales and thereafter combines those solutions to solve at increased scales. However, there were no multiplications nor complex \emph{twiddle} factors involved and the approach relied exclusively on integer arithmetics to achieve its goal. By working at multiple scales, and owing to the symmetry of the problem, intermediate computations can be reused, and thus, preventing the double computation of any sum and reducing the computational load from O($N^3$) to O($N^2 \log N$).

To enable this, the key is to define a loose discrete line that traverses the domain by visiting only integer positions, and therefore, not traversing in a straight way. Ascensions are defined recursively, making lines comprising two halves, which in turn separately come from two other lines of half their size and so on until lines that join only two points are reached, and the problem cannot be further reduced.

The authors of DRT eluded any sort of trigonometric relationship between the $x$ and $y$ variables, instead decomposing $u$ and $s$ variables of eq. \eqref{eq:radon2D} in binary and mixing them simultaneously (an index of both $u$ and $s$) at the binary level. This method also avoids the use of multiplications, which would have produced decimal numbers, as in the following: 
\begin{equation}
\stackrel{\scriptscriptstyle \sim}{f}\!(s,\;d) = \sum_{\mathbf{u} \in \mathbb{Z}^n_2} f(\lambda(\mathbf{u}), \; \; l^n_s(\mathbf{u})+d)
\label{eq:GotzfDRT} 
\end{equation}
This is the discrete version of eq. \eqref{eq:radon2D}, with the definition of discrete line $l^n_s(\mathbf{u})$ still pending. The auxiliary function, $\lambda(u_0, \ldots, u_{n-1}) = \sum_{i=0}^{n-1} u_i \cdot 2^i$, is used to convert binary multidimensional indexes to a decimal unidimensional index. $n$ is now $\log_2(N)$ and $N$ is the length of a side of our discrete domain.

In the work conducted by Marichal \emph{et al.}\cite{Marichal}, the formulation of DRT was modified to enable an extension to more dimensions. This study will adhere to that notation. The details on the formulation of the DRT can be found there. The following are the equations governing the DRT, the definition of discrete lines: eq. (\ref{eq:digitalLine}); the definition of a partial transform until stage $m$: eq. (\ref{eq:our2Dfm}); and the mapping between two stages: eq. (\ref{eq:map2D}):
\begin{multline}
\label{eq:digitalLine}
l^n_s(u_0, \ldots, u_{n-1}) = l^{n-1}_{\lfloor s/2 \rfloor}(u_0, \ldots, u_{n-2}) + 
u_{n-1} \left\lfloor \frac{s+1}{2} \right\rfloor =  \sum_{i=0}^{n-1} u_{n-1-i} \cdot \left\lfloor \frac{\frac{s}{2^i}+1}{2} \right\rfloor  
\end{multline}
\begin{multline}
\label{eq:our2Dfm}
\stackrel{\scriptscriptstyle \sim}{f}\!{}^m ( \overbrace{\, s_{n-m}, \overbrace{\, s_{n-m+1},\ldots,s_{n-1}\, }^{\boldsymbol{\sigma}}}^{\boldsymbol{s}}
| \;\overbrace{\, v_{m}, \ldots, v_{n-1} \, }^{\mathbf{v}} \, | \;d) = 
\sum_{\mathbf{u} \in \mathbb{Z}^{m}_2} 
f(\lambda(\mathbf{u}, \mathbf{v}) | l^m_{\lambda(\boldsymbol{s})}(\mathbf{u}) + d)
\end{multline}
\begin{multline}
\label{eq:map2D}
\stackrel{\scriptscriptstyle \sim}{f}\!{}^{m{+}1}
( \underbrace{\overbrace{\, s_{n-m-1}, \overbrace{\, s_{n-m},\ldots,s_{n-1}\,}^{\boldsymbol{\sigma}: \, m \; bits}}^{\boldsymbol{s}: \, m{+}1 \; bits}
| \overbrace{\; v_{m+1}, \ldots, v_{n-1} \,}^{\mathbf{v}: \, n{-}m{-}1\; bits} }_{n \; bits} | \,d) = \\
\stackrel{\scriptscriptstyle \sim}{f}\!{}^{m} (\boldsymbol{\sigma} | \, 0, \mathbf{v} | \, d) \; + \;
\stackrel{\scriptscriptstyle \sim}{f}\!{}^{m} (\boldsymbol{\sigma} | \,1, \mathbf{v} | \, d + s_{n-m-1} + \lambda(\boldsymbol{\sigma}))
\end{multline}

Notably, a single comma (,) is used to separate binary indexes, and a vertical bar ($|$) is used to separate different parameters.

Additionally, the number of bits in partial stages varies depending on $m$, the current stage. When $m = 0$, the array $\stackrel{\scriptscriptstyle \sim}{f}\!{}^{0}(s|v|d)$ is bidimensional, as variable $s$ is still empty. Thus, $\stackrel{\scriptscriptstyle \sim}{f}\!{}^{0}(-|v|d)$ maps directly to $f(x | y)$. When $m=n$, the last stage, variable $v$ will be emptied, and therefore, $\stackrel{\scriptscriptstyle \sim}{f}\!{}^{n}(s | - | d)$ is the desired result $\mathfrak{R}f(s | d)$. That can also be confirmed by evaluating the definition of the partial transform in stage $n$:  $\mathfrak{R}f(s | d) = \stackrel{\scriptscriptstyle \sim}{f}\!{}^n ( \, s_0, \ldots, s_{n-1}\, | \,d) = \sum_{u=0}^{N-1} 
f(\,u \,| \, l^n_{\lambda(\boldsymbol{s})}(u) + d).$ This last equation is the discrete version of the Radon transform as expressed in eq. (\ref{eq:GotzfDRT}) with multiplication of the slope substituted by the discrete line interaction between $u$ and $s$.

\subsection{Full sinogram construction}

Equations (\ref{eq:digitalLine}), (\ref{eq:our2Dfm}), and (\ref{eq:map2D}) comprise the core of the DRT algorithm for the basic quadrant with angles between 0 and $45^\circ$. A discrete version of the full sinogram can be accomplished, while reusing the algorithm that solves a quadrant, by applying it three times to versions of the input where axes are swapped and/or flipped conveniently, and their partial outputs are merged in a four times bigger global result. A typical output was shown at the bottom of figure \ref{fig:ConventionalShears}.

\section{SUMMATION OF PLANES IN 3D}
\label{sec:planes3D}

The equations for a 3D Radon transform in the continuum case are expressed as
\begin{equation}
    \mathfrak{R}f(\Pi) = \int_\Pi f(\boldsymbol{u}) d\boldsymbol{u} \label{eq:3DDRTsingle}
\end{equation}
where $\Pi$ will parameterize the plane to be considered. As we are adhering to the slope-intercept notation with $slope=\frac{\lvert s\rvert}{N-1} \leq 1$, we must split the formulae into three cases, depending on the axis near which the plane normals will circulate, as
\begin{multline}
 \mathfrak{R}_z f(slope_x, \,slope_y, \,d) = 
 \!{}\!{}\int \!{}\!{}\int \!{}\!{}f(u_x,\; u_y, \; u_x \cdot slope_x \: +  \: u_y\cdot slope_y \:+ \: d) \; \mathrm{d} u_x \, \mathrm{d} u_y, \\
 \textrm{or} \quad \mathfrak{R}_y f(slope_x, \,slope_z, \,d) = 
\!{}\!{}\int \!{}\!{}\int \!{}\!{}f(u_x, \; u_x \cdot slope_x \: + \: u_z\cdot slope_z  \:+ \: d,\; u_z) \; \mathrm{d} u_x \, \mathrm{d} u_z, \\
\label{eq:3Dcontinuous}
  \textrm{or} \quad \mathfrak{R}_x f(slope_y, \,slope_z, \,d) = 
 \!{}\!{}\int \!{}\!{}\int \!{}\!{}f(u_y \cdot slope_y \: + \: u_z\cdot slope_z \:+ \: d,\;u_y,\; u_z) \; \mathrm{d} u_y \, \mathrm{d} u_z.
\end{multline}
with the first case considering the planes of form $z- x\;slope_x - y \; slope_y - d = 0$, and so on. For simplicity, we will expose the basic algorithm just for this case.

\subsection{3D multiscale DRT of planes}
The discrete approximation that we will compute is as follows:
\begin{multline}
\stackrel{\scriptscriptstyle \sim}{f}\!(s_1 \, | \, s_2 \, |\,d) = 
 \sum_{\mathbf{u_1} \in \mathbb{Z}^n_2}
 \sum_{\mathbf{u_2} \in \mathbb{Z}^n_2}  f(\lambda(\mathbf{u_1})| \, \lambda(\mathbf{u_2}) | \, l^n_{s_1}(\mathbf{u_1}) + l^n_{s_2}(\mathbf{u_2}) +d).
\label{eq:3DfDRT}
\end{multline}

If we define the solution up to stage $m$ as
\begin{multline}
\label{eq:3Dfm}
\stackrel{\scriptscriptstyle \sim}{f}\!{}^m ( \overbrace{\, s_{1_{n-m}}, \ldots,s_{1_{n-1}}\, }^{\mathbf{s_1}}
| \;\overbrace{\, v_{1_{m}}, \ldots, v_{1_{n-1}} \, }^{\mathbf{v_1}} \, | 
\overbrace{\, s_{2_{n-m}}, \ldots,s_{2_{n-1}}\, }^{\mathbf{s_2}}
| \;\overbrace{\, v_{2_{m}}, \ldots, v_{2_{n-1}} \, }^{\mathbf{v_2}} | \;d) = \\
\sum_{\mathbf{u_1} \in \mathbb{Z}^{m}_2} \sum_{\mathbf{u_2} \in \mathbb{Z}^{m}_2} 
f(\lambda(\mathbf{u_1}, \mathbf{v_1}) \, | \; 
  \lambda(\mathbf{u_2}, \mathbf{v_2}) \, | 
 l^m_{\lambda(\mathbf{s_1})}(\mathbf{u_1}) + l^m_{\lambda(\mathbf{s_2})}(\mathbf{u_2}) + d), 
\end{multline}
the solution is accomplished by applying the following mapping recursively:
\begin{multline}
\label{eq:map3D}
\stackrel{\scriptscriptstyle \sim}{f}\!{}^{m{+}1}
(  s_{1_{n-m-1}}, \overbrace{\, s_{1_{n-m}},\ldots,s_{1_{n-1}}\,}^{\boldsymbol{\sigma_1}} |\;
  \overbrace{\; v_{1_{m+1}}, \ldots, v_{1_{n-1}} \,}^{\mathbf{v_1}} |  \\
   s_{2_{n-m-1}}, \overbrace{\, s_{2_{n-m}},\ldots,s_{2_{n-1}}\,}^{\boldsymbol{\sigma_2}} | \;
  \overbrace{\; v_{2_{m+1}}, \ldots, v_{2_{n-1}} \,}^{\mathbf{v_2}} | \; d) = \\
  {\color{blue}
  \sum_{v_{1_m}}^{0,1} \sum_{v_{2_m}}^{0,1} \stackrel{\scriptscriptstyle \sim}{f}\!{}^{m} \Bigl(\boldsymbol{\sigma_1} \, | \, v_{1_m}, \mathbf{v_1} \, | \, \boldsymbol{\sigma_2} \, | \, v_{1_m}, \mathbf{v_2} \, |\,  d + v_{1_m} \cdot (s_{1_{n-m-1}} + \lambda(\boldsymbol{\sigma_1})) + v_{2_m} \cdot (s_{2_{n-m-1}} + \lambda(\boldsymbol{\sigma_2}))\Bigr) = } \\
\stackrel{\scriptscriptstyle \sim}{f}\!{}^{m} (\boldsymbol{\sigma_1} \, | \, 0, \mathbf{v_1} \, | \, \boldsymbol{\sigma_2} \, | \, 0, \mathbf{v_2} \, |\,  d) + 
\stackrel{\scriptscriptstyle \sim}{f}\!{}^{m} (\boldsymbol{\sigma_1} \, | \, 1, \mathbf{v_1} \, | \, \boldsymbol{\sigma_2} \, | \, 0, \mathbf{v_2}) \, | \, d + s_{1_{n-m-1}} + \lambda(\boldsymbol{\sigma_1})) + \\
\stackrel{\scriptscriptstyle \sim}{f}\!{}^{m} (\boldsymbol{\sigma_1} \, | \, 0, \mathbf{v_1} \, | \, \boldsymbol{\sigma_2} \, | \, 1, \mathbf{v_2} \, | \, d + s_{2_{n-m-1}} + \lambda(\boldsymbol{\sigma_2})) + \\
\stackrel{\scriptscriptstyle \sim}{f}\!{}^{m} (\boldsymbol{\sigma_1} | 1, \mathbf{v_1} | \boldsymbol{\sigma_2} | \, 1, \mathbf{v_2} | 
\, d + s_{1_{n-m-1}} + \lambda(\boldsymbol{\sigma_1}) + s_{2_{n-m-1}} + \lambda(\boldsymbol{\sigma_2})).
\end{multline}
As stated previously, these formulae were already derived\cite{Marichal}, and translate directly into an algorithm, which is given in Algorithm \ref{alg:3DDRT}. However it only solves a reduced number of planes: those in the form  $z- x\;slope_x - y \; slope_y - d = 0$, that is, a group of planes whose normals are $\mathbf{n} = (-\frac{s_1}{N-1}, -\frac{s_2}{N-1}, 1)$ with $\frac{s_1}{N-1}$ and $\frac{s_2}{N-1} \in [0, \,1]$. 

\begin{algorithm}
\caption{Compute the 3D DRT of a dodecant}
\label{alg:3DDRT}
\textbf{Input:} $f(x, y, z)$ consisting of $N\times N \times N$ data\\
\textbf{Output:} 3D discrete Radon transform of f, $\mathfrak{R}f(s1, s2, d)$ consisting of $N\times N \times 3N$ data\\
\textbf{function} basicDRT3D(f)
\begin{algorithmic}
\State $n \gets \mathrm{\mathbf{log2}}(N)$
\State ${f^m} \leftarrow \mathrm{\mathbf{zeros}}(N, \: N, \: 3\cdot N)$
\State ${f^m}(0:N-1, \: 0:N-1, \: 2N:3N-1) \leftarrow {f}(0:N-1, \: 0:N-1, \: 0:N-1)$
\State ${f^{m+1}} \leftarrow \mathrm{\mathbf{zeros}}(N, \: N, \: 3\cdot N)$
\For{$m=0$ to $n-1$}
	\For{$d=0$ to $3\cdot N-1$}
    	\For{$v_1=0$ to $(1 << (n-m-1) )-1$}
			\For{$\sigma_1=0$ to $(1 << m)-1$}
			\For{$v_2=0$ to $(1 << (n-m-1) )-1$}
			\For{$\sigma_2=0$ to $(1 << m)-1$}
			\For{$s_{1-lsb} = 0, 1$}
			\For{$s_{2-lsb} = 0, 1$}
            	\State $f_{00} \leftarrow f^m(\sigma_1 + (v_1 << (m+1)),\: \ldots$
            	\State $\ldots \sigma_2 + (v_2 << (m+1)),\: d )$
            	\LineComment{Next 3 accesses should be 0 ...}
            	\LineComment{...if overflow 3rd dimension}
            	\State $f_{10} \leftarrow f^m(\sigma_1 + (1 << m) + (v_1 << (m+1)),\ldots$
            	\State $\ldots \: \sigma_2 + (v_2 << (m+1)),\: d + \sigma_1 + s_{1-lsb})$ 
            	\State $f_{01} \leftarrow f^m(\sigma_1 + (v_1 << (m+1)), \ldots$
            	\State $\ldots \: \sigma_2 + (1 << m) + (v_2 << (m+1)), \ldots$
            	\State $ \ldots \: d + \sigma_2 + s_{2-lsb})$
            	\State $f_{11} \leftarrow f^m(\sigma_1 + (1 << m) + (v_1 << (m+1)), \ldots$
            	\State $\quad \ldots \sigma_2 + (1 << m) + (v_2 << (m+1)), \ldots$
            	\State $\ldots \: d + \sigma_1 + \sigma_2 + s_{1-lsb} + s_{2-lsb} )$
            	\LineComment{After 4 reads, we can write 1 datum}
            	\State $f^{m+1}(s_{1-lsb} + (\sigma << 1) + (v << (m+1)), \ldots$
            	\State $\ldots \: s_{2-lsb} + (\sigma_2 << 1) + (v_2 << (m+1)), \: d) \ldots$ 
            	\State $\qquad \quad \ldots \leftarrow f_{00}+f_{10}+f_{01}+f_{11}$ 
            	\EndFor \Comment{variable $s_{2-lsb}$}
            	\EndFor \Comment{variable $s_{1-lsb}$}
                \EndFor \Comment{variable $\sigma_2$}
                \EndFor \Comment{variable $v_2$}
            \EndFor \Comment{variable $\sigma_1$}
        \EndFor \Comment{variable $v_2$}
    \EndFor \Comment{variable $d$}
    \State ${f^m} \leftarrow {f^{m+1}}$ \Comment{Reuse buffers interchanging last output for next input}
\EndFor \Comment{variable $m$}
\State $\mathfrak{R}f  \leftarrow f^m$ \\
\Return{${\mathfrak{R}f}$}
\end{algorithmic}
\end{algorithm}

\begin{figure}[h]
\centering
\includegraphics[width=.35\textwidth]{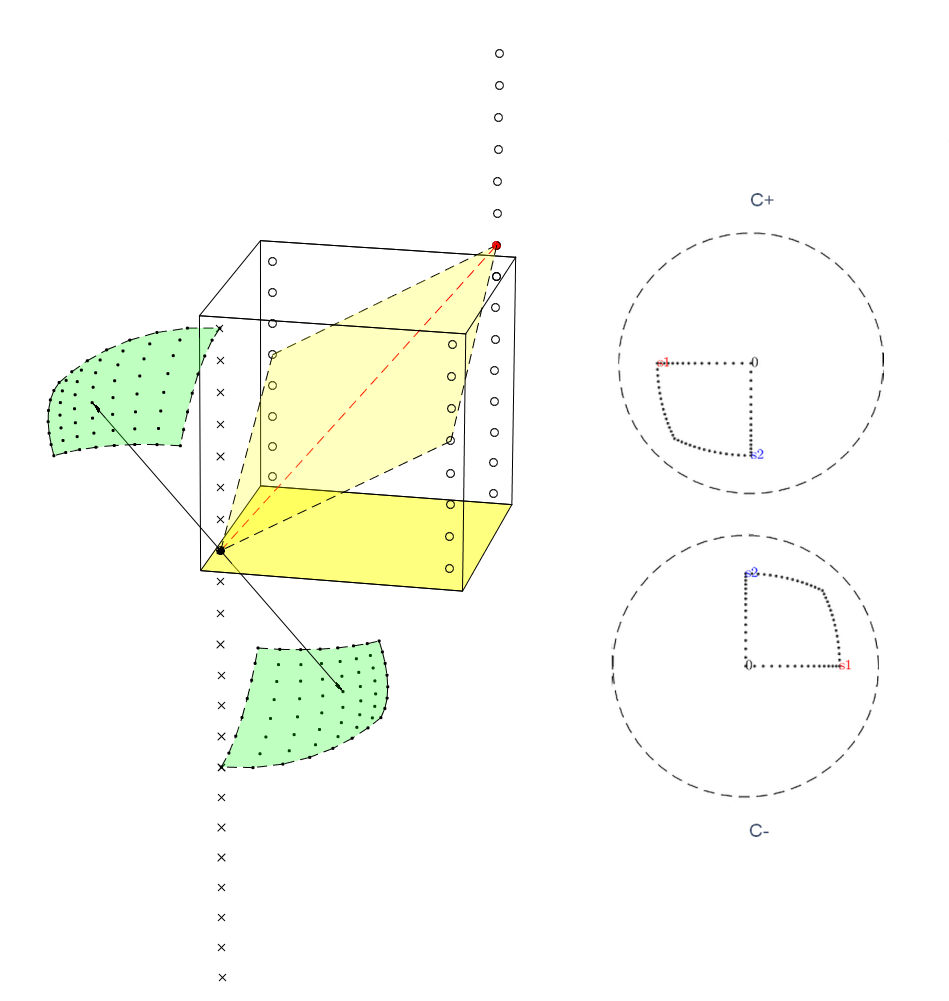} 
\caption{Right: Depiction of the parameters of the basic algorithm of 3D DRT. Left: Boundary x-y coordinates of plane normals covered by the basic algorithm, as projected into positive and negative hemispheres cut through $z=0$.}
\label{fig:1dodecante}
\end{figure}

Figure \ref{fig:1dodecante} depicts the planes, normals, slopes, and displacements considered by the basic algorithm for 3D DRT.
The depiction symbolizes the algorithm for $\mathbf{p_0}=(0,0,0)$ and ascents in direction $\mathbf{d} = (0,0,1)$, that is, pointing to the positive $z-$axis. Those are the \emph{origin point} and \emph{null normal} emerging from it. The discrete computed planes are defined by three points, $\mathbf{p_0},\: \mathbf{p_1}+s_1\cdot\mathbf{d}$ and $\mathbf{p_2}+s_2\cdot\mathbf{d}$; with $\mathbf{p_1}=(1, 0, 0)\cdot (N-1)$ and $\mathbf{p_2}=(0, 1, 0)\cdot (N-1)$, which are adjacent neighbors of $\mathbf{p_0}$ in the face $z=0$ of the cube of size $N \times N\times N$. 

{\color{blue}
In the figure, the origin point $\mathbf{p_0}$ is represented as a black dot, whereas the null plane is represented by the almost opaque yellow polygon; the cube face $z=0$. When the points $\mathbf{p_1}+s_1\cdot\mathbf{d}$ and $\mathbf{p_2}+s_2\cdot\mathbf{d}$ ascend, the normal through $\mathbf{p_0}$ deviates from the null normal and the planes to be considered differ as well from the null plane. Without moving $\mathbf{p_0}$, there are $N \times N$ planes to be computed; those arising as $s_1$ and $s_2$ vary from $0$ to $N-1$. The case for $s_1 = s_2 = 4$ is represented as the semitransparent yellow plane. The arrows through $\mathbf{p_0}$ represent the normal in that case, and the other possibilities are represented as the dots drawn on the green spherical caps; there are two of them, indicating that normals can be assumed to be positive or negative.

The crosses that form a column moving up and down from $\mathbf{p_0}$ in the direction of $\mathbf{d}$ are the other different displacements of $\mathbf{p_0}$ to be considered; additional $N-1$ positions upwards and $2(N-1)$ downwards. The columns of circles in the other vertices of the null face, however, are the positions that can be reached without moving $\mathbf{p_0}$ from its null displacement. All this combinations are already considered by the \emph{for} loops in the Algorithm \ref{alg:3DDRT}, and they explain the size of the output: there will be two dimensions, each of size $N$ for describing the normals; along with one dimension for the displacements of size $3N-2$, which for the sake of memory allocation simplicity can be assumed to be $3N$ in the algorithm description.}


This basic algorithm covers a \emph{dodecant}, a twelfth part of a hemisphere of normals. 
\subsection{Shape and construction of the output}
\begin{figure}[h]
\centering
\includegraphics[width=.6\textwidth]{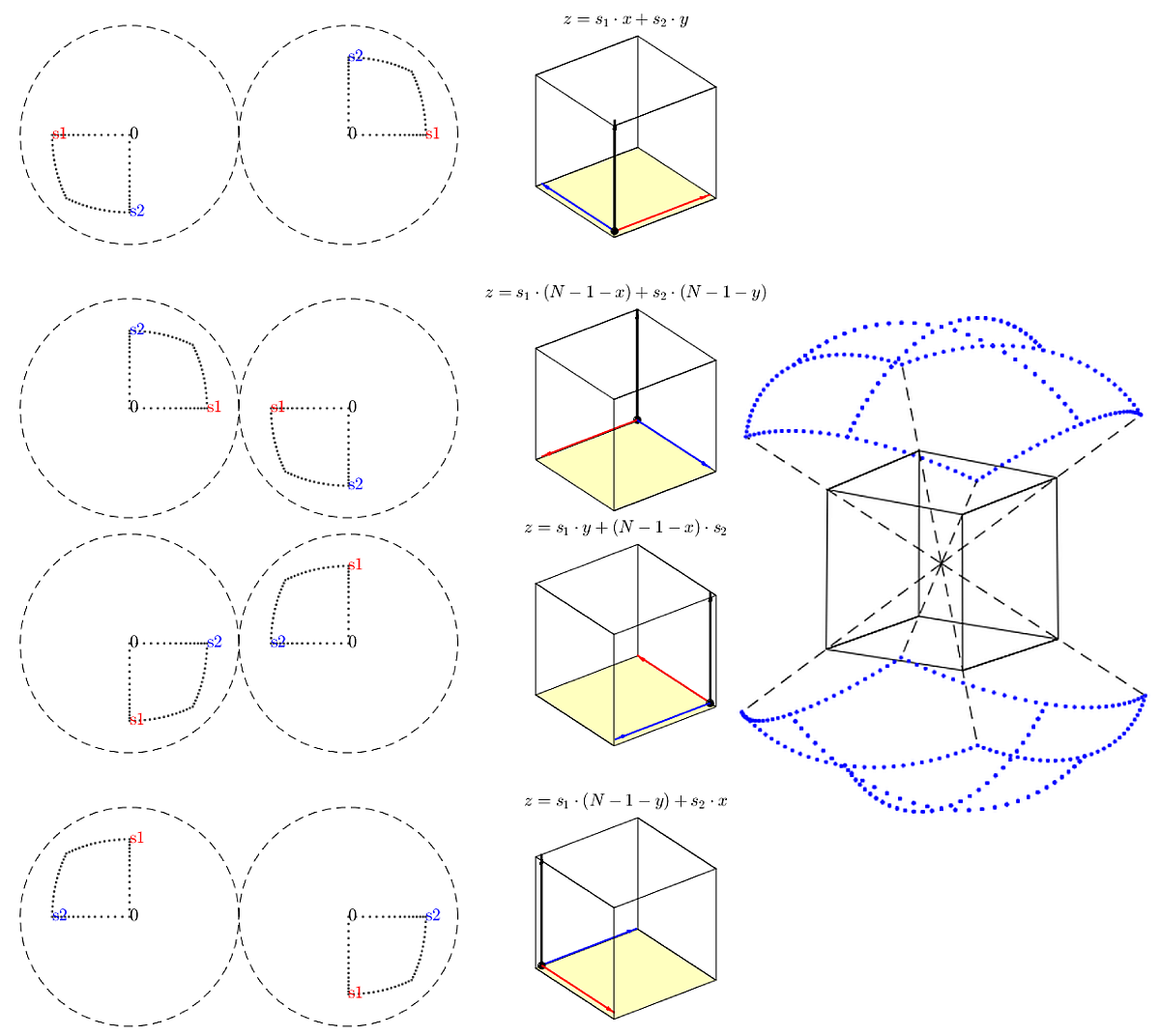} 
\caption{Depiction of the four dodecants that compute the entire set of planes near the $z-$normal: null plane (face in yellow), origin point $\mathbf{p_0}$ (black circle), direction of ascents (black arrow), direction of neighbors points $\mathbf{p_1}$ and $\mathbf{p_2}$ (red and blue arrows), and location of boundary $x-y$ coordinates for plane normals (dots) within positive and negative $z-$hemispheres (dashed circles).
Left: Quadrilateralized spherical caps emerge from the z-faces, inspired from Chan \& O'Neill (1975).}
\label{fig:az_el_z}
\end{figure}

To compute the rest of the hemisphere, eleven more cases remain to be considered. Three similar versions of the basic algorithm are needed to cover the possibilities near null normal in direction $\mathbf{d} = (0,0,1)$. These additional cases correspond to placing the origin $\mathbf{p_0}$ on the remainder of vertices on the face cube at $z=0$: $(0, 1, 0), (1, 0, 0)$, and $(1, 1, 0)$, scaled by $(N-1)$. In each case, $\mathbf{p_1}$ and $\mathbf{p_2}$ are the face cube vertices adjacent to $\mathbf{p_0}$ such that the right-hand rule is verified: $(\mathbf{p_1}-\mathbf{p_0}) \times (\mathbf{p_2}-\mathbf{p_0}) = \mathbf{d}$ (see figure \ref{fig:az_el_z}).
The four runs will achieve all the planes whose normals are in the spherical patch corresponding to face cube $z=0$ of the quadrilateralized spherical cube \cite{ChanOneill} (see far right of figure \ref{fig:az_el_z}). We will refer to them as the $z-$normal set of planes.

\begin{figure}[h]
\centering
\includegraphics[width=\textwidth]{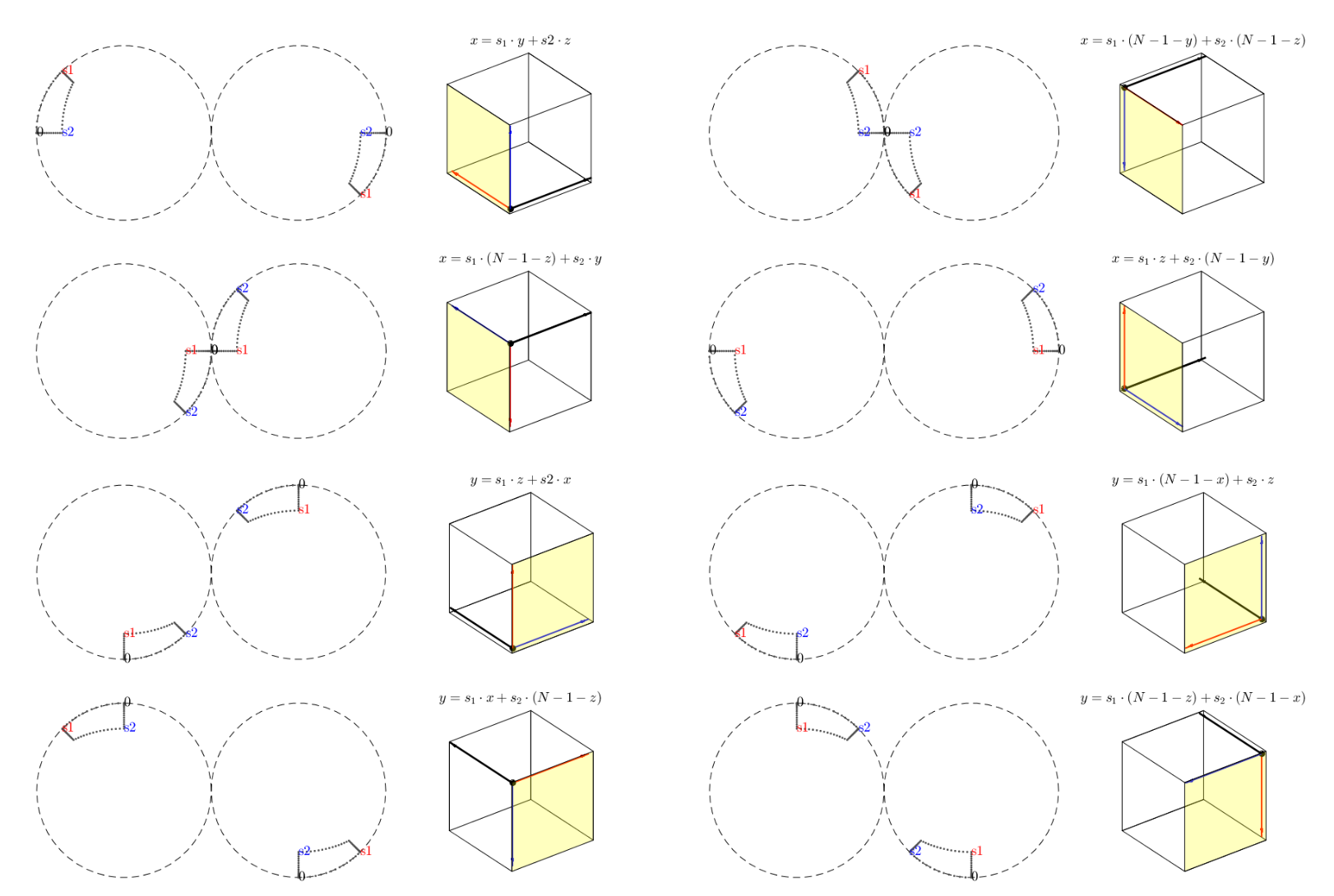} 
\caption{Depiction of the $x-$normal and $y-$normal set of planes, equivalent to figure \ref{fig:az_el_z} for the $z-$normal set.}
\label{fig:az_el_xy}
\end{figure}

The computation of the $y-$normal and $x-$normal set of planes, discrete versions of eqs. \eqref{eq:3Dcontinuous}, will require eight additional runs of variations of Algorithm \ref{alg:3DDRT}. As with 2D DRT, instead of modifying the algorithm, it is easier to operate permuting and/or flipping dimensions on the input, as dictated by the formula above each cube in figure \ref{fig:az_el_xy}.

The permutations and flipping of dimensions of the input to achieve the computation of all the dodecants are made explicit in Algorithm \ref{alg:3DDRT_cube}. Note that four dodecants inversely traverse $\mathbf{d}$ with respect to the figures because of symmetry.

\begin{algorithm}
\caption{Compute the 3D DRT of a cube}
\label{alg:3DDRT_cube}
\textbf{Input:} $f(x, y, z)$ consisting of $N\times N \times N$ data\\
\textbf{Output:} $\mathfrak{R}f$ consisting of $4N \times 4N \times 3N-2$ data\\
\textbf{function} global3DDRT(f)
\begin{algorithmic}
\State $InOrder \gets [1, 2, 3; -2, 1, 3; -1, -2, 3; 2, -1, 3;$ \Comment{In these arrays initializations}
\State $\quad \qquad \qquad -1, -3, -2; -3, -1, 2; 3, -1, -2; 1, 3, -2;$ \Comment{(,) separates columns}
\State $\quad \qquad \qquad -3, -2, -1; -2, -3, 1; -2, 3, -1; 3, 2, -1]$ \Comment{(;) separates rows}
\State $OutOrder \gets [-2, -1, 3; -1, 2, 3; 2, 1, 3; 1, -2, 3;$
\State $\quad \qquad \qquad 2, 1, 3; 1, -2, 3; -1, 2, 3; -2, -1, 3;$
\State $\quad \qquad \qquad 2, 1, 3; -1, 2, 3; 1, -2, 3; -2, -1, 3]$
\State $DodecantCoords \gets [1, 1; 1, 2; 2, 2; 2, 1;$
\State $\qquad \qquad \qquad \qquad \quad 0, 2; 0, 1; 3, 2; 3, 1;$
\State $\qquad \qquad \qquad \qquad \quad 2, 0; 1, 0; 2, 3; 1, 3]$
\State $N \gets \mathrm{size}(f, 1)$ \Comment{Extract size from any dimension of input}
\State $\mathfrak{R}f \gets \mathrm{zeros}(4N, 4N, 3N-2)$ \Comment{Memory allocation for the output}
\State \For{$k$ = 0 to 11}  \Comment{Compute and merge the twelve dodecants}
\State $inOrder \gets InOrder(k, \,0:2)$
\State $f \gets \mathrm{permute(}f, \textrm{abs(}inOrder))$ \Comment{Permute input dimensions}
\For{dim = 0 to 2}
    \If{$inOrder(dim) < 0$}
    \State $f \gets \mathrm{flip(}f, \, dim)$ \Comment{Flip input dimensions}
    \EndIf
\EndFor
\State \LineComment{Perform the computation of a dodecant}
\State $RfDodecant \gets \mathrm{basicDRT3D}(f)$
\State \LineComment{Prepare and merge the dodecant output into global solution}
\State $outOrder \gets [OutOrder(k, \, 0:2)]$ 
\State $RfDodecant \gets \mathrm{permute(}RfDodecant, \, \textrm{abs(}outOrder))$ \Comment{Permute}
\For{$dim \: = \: 0 \: to \: 2$}
    \If{$outOrder(dim) < 0$}
    \State $RfDodecant \gets \mathrm{flip}(RfDodecant, \, dim)$ \Comment{Flip output}
    \EndIf
\EndFor
\State \LineComment{Place RfDodecant in $\mathfrak{R}f$ according to dodecantCoords}
\State $dC \gets DodecantCoords(k, \,0:1)$
\State $\mathfrak{R}f(dC(0)\cdot N:(dC(0)+1)\cdot N-1, \: dC(1)\cdot N:(dC(1)+1)\cdot N-1, \: 0:3N-2) \gets RfDodecant$
\EndFor  \\
\Return{${\mathfrak{R}f}$}
\end{algorithmic}
\end{algorithm}

The output of the global algorithm will occupy $4\cdot N \times 4\cdot N \times 3N-2$ in memory. The particular shape of a dodecant is bounded by $N\times N\times 3N-2$. But several of these data will always be zero as there will be displacements from where a certain $s_1, \, s_2$ never touches the $N \times N \times N$ portion populated with data, in particular, non-null values are within the $(0:N-1, 0:N-1, 0:N-1+s_1+s_2)$.
 
In addition, despite reserving $4\times 4$ space for dodecants, only 12 of them contain data. 

\begin{figure}[h]
\centering
\includegraphics[width=.7\textwidth]{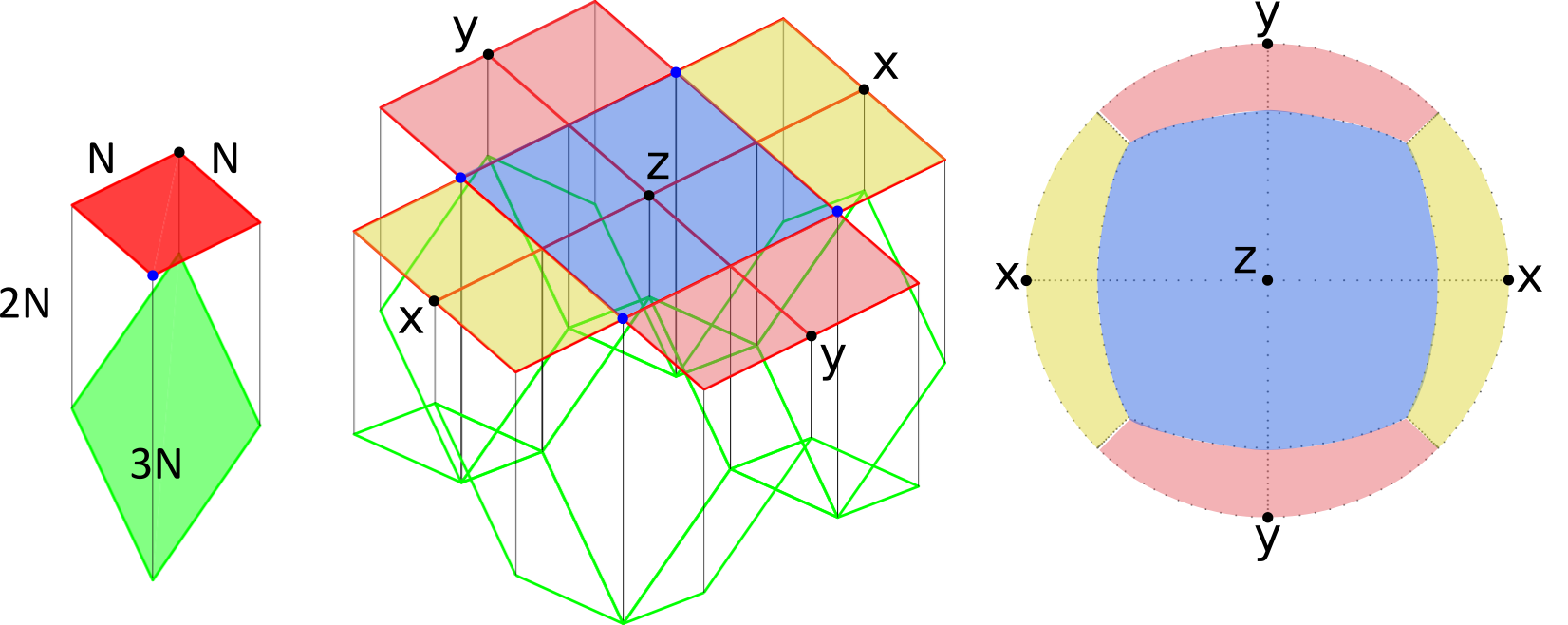}
\caption{\color{blue} Shape of the output. Left: Isolated dodecant output is shown. The black dot is located on top of the column with $s_1 = 0$ and $s_2 = 0$, which will hold at most $N$ non-null displacements; the blue dot is located on top of the column with $s_1 = N-1$ and $s_2 = N-1$, which can hold $3N-2$ non-null displacements.
Middle: 12 dodecants are merged. Each individual dodecant output is rotated so that all of them fit accurately. The black dots mark the position of $x, y$, and $z$ null normals, accordingly to hemisphere of normals, depicted at right.}
\label{fig:shapeOutput}
\end{figure}

The shape of an individual dodecant and global output is shown in figure \ref{fig:shapeOutput}. An actual run of the algorithm for a cube of size $64^3$ with non-null entries only in the cube corners and the center is shown in figure \ref{fig:exampleOutput}.

\section{SUMMATION OF LINES IN 3D}
\label{sec:lines3D}
The equations for the 3D John transform in the continuum case is expressed as
\begin{equation}
    \mathfrak{J}f(\Lambda) = \int_\Lambda f(\boldsymbol{u}) d\boldsymbol{u} \label{eq:3DDJTsingle}
\end{equation}
where $\Lambda$ will parameterize the line to be considered. We will adhere to the slope-intercept notation, with $slope=\frac{\lvert s\rvert}{N-1} \leq 1$. The previous formula leads to three cases, depending on the axis with a lesser angle with respect to the line vector:
\begin{multline*}
\mathfrak{J}_x f(s_y, s_z, d_y, d_z) = \!{}\!{}\int \!{}\!{}f(u_x,\; \,s_y\cdot u_x + d_y, \,s_z \cdot u_x + d_z) \; \mathrm{d} u_x, \\
\mathfrak{J}_y f(s_x, s_z, d_x, d_z) = \!{}\!{}\int \!{}\!{}f(s_x\cdot u_y + d_x, \, u_y, \,s_z \cdot u_y + d_z)  \mathrm{d} u_y, \\
\mathfrak{J}_z f(s_x, s_y, d_x, d_y) = \!{}\!{}\int \!{}\!{}f(s_x\cdot u_z + d_z, \,s_y \cdot u_z + d_y, \; u_z) \; \mathrm{d} u_z.
\end{multline*}

We can derive the algorithm just for the first case, in which the lines adopt the parametric form \[\begin{cases}
        x=u\\
        y=slope_y\cdot u+ d_y\\
        z=slope_z\cdot u+ d_z  \qquad \qquad \textrm{or in vector form}\\
    \end{cases}\] \[(x,y,z)=(0,\, d_y, \, d_z) + u \cdot (1, slope_y, slope_z).\]

\subsection{3D multiscale DJT of lines}
The discrete approximation that we will compute is
\begin{equation}
\stackrel{\scriptscriptstyle \sim}{f}\!(s_1 | \; s_2 | \; d_1 |\;d_2) = \\
 \sum_{\mathbf{u} \in \mathbb{Z}^n_2}
f(\lambda(\mathbf{u})| \; l^n_{s_1}(\mathbf{u}) + d_1 | \; l^n_{s_2}(\mathbf{u}) +d_2).
\label{eq:3DfDJT}
\end{equation}

If we define the solution up to stage $m$ as
\begin{multline}
\label{eq:3DJfm}
\stackrel{\scriptscriptstyle \sim}{f}\!{}^m ( \overbrace{\, v_{m}, \ldots, v_{n-1} \, }^{\mathbf{v}} \, |
\overbrace{\, s_{1_{n-m}}, \ldots,s_{1_{n-1}}\, }^{\mathbf{s_1}} | \; 
\overbrace{\, s_{2_{n-m}}, \ldots,s_{2_{n-1}}\, }^{\mathbf{s_2}} | \;
d_1 | \; d_2) = \\
\sum_{\mathbf{u} \in \mathbb{Z}^{m}_2} 
f(\lambda(\mathbf{u}, \mathbf{v}) \, | \;
 l^m_{\lambda(\mathbf{s_1})}(\mathbf{u}) + d_1 | \;
 l^m_{\lambda(\mathbf{s_2})}(\mathbf{u}) + d_2), 
\end{multline}
the solution is obtained by applying the following mapping recursively.
\begin{multline}
\label{eq:map3DJ}
\stackrel{\scriptscriptstyle \sim}{f}\!{}^{m{+}1}
(\overbrace{\; v_{m+1}, \ldots, v_{n-1} \,}^{\mathbf{v}} |
s_{1_{n-m-1}}, \overbrace{\, s_{1_{n-m}},\ldots,s_{1_{n-1}}\,}^{\boldsymbol{\sigma_1}} | \; 
s_{2_{n-m-1}}, \overbrace{\, s_{2_{n-m}},\ldots,s_{2_{n-1}}\,}^{\boldsymbol{\sigma_2}} |\; d_1 |\; d_2) = \\
\stackrel{\scriptscriptstyle \sim}{f}\!{}^{m} (0, \, \mathbf{v} \,| \, \boldsymbol{\sigma_1} \, | \, \boldsymbol{\sigma_2} \, |\, d_1 \, | \, d_2) + 
\stackrel{\scriptscriptstyle \sim}{f}\!{}^{m} (1, \, \mathbf{v} | \, \boldsymbol{\sigma_1} \, | \, \boldsymbol{\sigma_2} \, | \, 
d_1 \,+\,s_{1_{n-m-1}}\,+\,\lambda(\boldsymbol{\sigma_1})\, | \, 
d_2 \,+\,s_{2_{n-m-1}}\,+\,\lambda(\boldsymbol{\sigma_2})) .
\end{multline}

As expected, and in accordance with figure \ref{fig:schemes}, there will be four possible combinations and each sums two values.

This can be considered to be a five-dimensional transform at partial stages, where $v$ is being emptied and its bits transferred and doubled at each stage, making $s_1$ and $s_2$ grow. Therefore, the transform starts becoming 3D, and $\stackrel{\scriptscriptstyle \sim}{f}\!{}^0 ( \overbrace{ \mathbf{v} }^{n \, \mathrm{bits}} \, | \overbrace{\, \cancel{\mathbf{s_1}} }^{0 \, \mathrm{bits}} | \, \overbrace{\cancel{\mathbf{s_2}}\, }^{0 \, \mathrm{bits}} | \; d_1 \,| \;d_2)$ maps directly to $f(x ,\,y,\,z)$. However, after $n$ stages, it becomes four-dimensional (4D), as desired: \[\stackrel{\scriptscriptstyle \sim}{f}\!{}^n ( \overbrace{ \cancel{\mathbf{v}} }^{0 \, \mathrm{bits}} \, | \overbrace{\, \mathbf{s_1} }^{n \, \mathrm{bits}} | \, \overbrace{\mathbf{s_2}\, }^{n \, \mathrm{bits}} | \; d_1 \,| \;d_2) = \stackrel{\scriptscriptstyle \sim}{f}\!(s_1| \;  s_2| \; d_1 |\;d_2).\]

Owing to simplicity and to be able to perform temporal buffers swapping after each stage, we will describe memory in the algorithm as 3D, with the least-significant dimension shared between $v, s_1$, and $s_2$, and thereafter completed with zeros to always occupy $2n$ bits: \[\stackrel{\scriptscriptstyle \sim}{f}\!{}^m ( \underbrace{\overbrace{ \mathbf{v} }^{n-m \, \mathrm{bits}}, \, \overbrace{\mathbf{s_1}}^{m \, \mathrm{bits}}, \, \overbrace{\mathbf{s_2}}^{m \, \mathrm{bits}} , \, \overbrace{\mathbf{0}}^{n-m \, \mathrm{bits}}   }_{2n \, \mathrm{bits}} | \, d_1 | \,d_2)\]
Notably, the memory indexes in a computer are normally interpreted as least significant to the right, but these formulae are interpreted with greater significance to the right to adhere to the original DRT authors' notation.

All these have been considered in Algorithm \ref{alg:3DDJT}. However, that only solves a twelfth of the total possible 3D lines cutting a cube. A similar discussion that was pertinent for plane normals in 3D DRT is now applicable to line vectors in 3D DJT. Other authors have named these as $x-$driven, $y-$driven, and $z-$driven sets of lines. The cube is then divided into three sectors, similar to reversed octahedrons comprising two square pyramids joined on their apexes, instead of their bases. We prefer to describe them as the projections of faces in the quadrilateralized spherical cube geometry. However, in the end, they represent a similar concept. To obtain the remainder of the dodecants, similar permutations and flipping of axes that were applied to the 3D DRT can now be applied, and are shown synthesized in the $InOrder$ array in Algorithm \ref{alg:3DDRT_cube}.

Merging the output of each dodecant into a global output is not easy or even useful because the output of this basic algorithm is now 4D. We will skip that part in this discussion. 

\begin{figure}
    \centering
    \includegraphics[width=.7\textwidth]{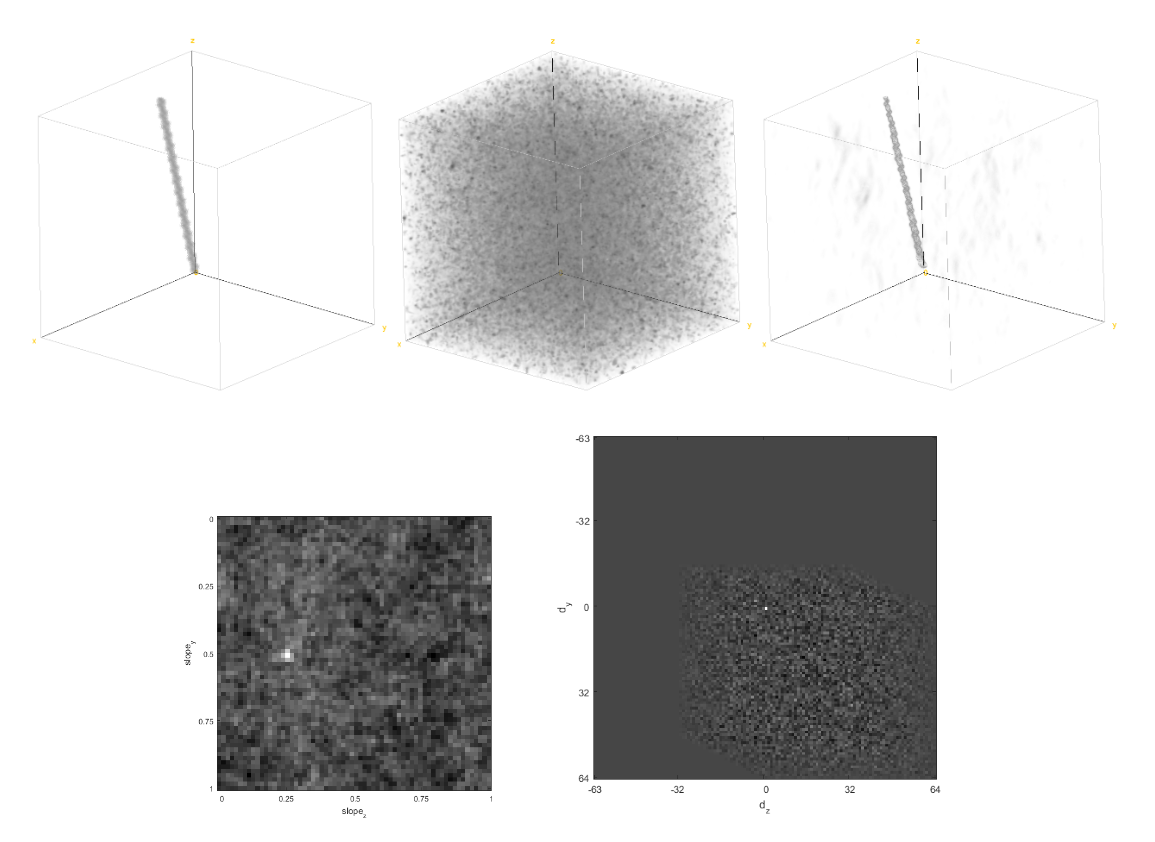}
    \caption{\color{blue} Top row, left to right: Voxels aligned in a 3D volume conforming a ray; these voxels become ``buried'' into white noise, and its 3D DJT is computed; the lower amplitude coefficients in the 3D DJT domain are suppressed, and then the signal is back-projected: the ray is not affected by the thresholding, while the noise disappears. Bottom: Two slices of the 4D output of the DJT transform before coefficient thresholding. Left: Displacements fixed to those of the ray. Right: Slopes fixed to those of the ray. In both cases, there is a magnitude peak in the ray parameters: those that will exceed the threshold.}
    \label{fig:output3DDJT}
\end{figure}

In figure \ref{fig:output3DDJT} we show an exemplary run of the basic dodecant algorithms. {\color{blue} We placed a ray in a cube and then added white noise that exceeds it in amplitude. However, the ray coefficient in John's domain, in its displacements-slopes coordinates, still stands out from its surroundings. Thus a thresholding of the highest amplitude coefficients, and subsequent backprojection, eliminates most of the noise and brings back the line visible again.}

\begin{algorithm}
\caption{Compute the 3D DJT of a dodecant}
\label{alg:3DDJT}
\textbf{Input:} $f(x, y, z)$ consisting of $N\times N \times N$ data\\
\textbf{Output:} 3D discrete John transform of $f$, $\mathfrak{J}f(s1, s2, d1, d2)$ consisting of $N\times N \times 2N \times 2N$ data\\
\textbf{function} basicDJT3D($f$)
\begin{algorithmic}
\State $n \gets \mathrm{\mathbf{log2}}(N)$
\State ${f^m} \leftarrow \mathrm{\mathbf{zeros}}(N^2, \: 2N, \: 2N)$
\State ${f^m}(0:N-1, \: N:2N-1, \: N:2N-1) \leftarrow {f}(0:N-1, \: 0:N-1, \: 0:N-1)$
\State ${f^{m+1}} \leftarrow \mathrm{\mathbf{zeros}}(N^2, \: 2N, \: 2N)$
\For{$m=0$ to $n-1$}
	\For{$d_2=0$ to $2N-1$}
	    \For{$d_1=0$ to $2N-1$}
	        \For{$\sigma_2=0$ to $(1 \ll m)-1$}
	            \For{$s_{2-lsb} = 0, 1$}
	                \For{$\sigma_1=0$ to $(1 \ll m)-1$}
    			        \For{$s_{1-lsb} = 0, 1$}
    	                    \For{$v=0$ to $(1 \ll (n-m-1) )-1$}
            	\State $f_{0} \leftarrow f^m((v \ll 1) + (\sigma_1 \ll n-m) \ldots$ 
            	\State $\ldots + (\sigma_2 \ll n), \; d_1, \; d_2 )$
            	\LineComment{Next access should be clamped to 0...}
            	\LineComment{...if exceeds 2nd or 3rd dimension limits}
            	\State $f_{1} \leftarrow f^m(1 + (v \ll 1) + (\sigma_1 \ll n-m) \ldots$
            	\State $\ldots + (\sigma_2 \ll n),\, d_1 \,+\, \sigma_1\, + \, s_{1-lsb}, \: \ldots$
            	\State $\ldots d_2 \,+\, \sigma_2 \, +\, s_{2-lsb} )$ 
            	\LineComment{After 2 reads, we can write 1 in next stage}
            	\State $f^{m+1}(v + (s_{1-lsb} \ll n-m-1) + \ldots$
            	\State $\ldots (\sigma_1 \ll n-m) + (s_{2-lsb} \ll n) + \ldots$
            	\State $\ldots (\sigma_2 \ll n+1),\: d_1 ,\:  d_2 ) = f_0 \,+\, f_1$ 
            	\EndFor \Comment{variable $v$}
            	\EndFor \Comment{variable $s_{1-lsb}$}
                \EndFor \Comment{variable $\sigma_1$}
                \EndFor \Comment{variable $s_{2-lsb}$}
            	\EndFor \Comment{variable $\sigma_2$}
                \EndFor \Comment{variable $d_1$}
                \EndFor \Comment{variable $d_2$}
    \State ${f^m} \leftarrow {f^{m+1}}$ \Comment{Reuse buffers interchanging last output for next input}
\EndFor \Comment{variable $m$}
\State $\mathfrak{J}f  \leftarrow \mathrm{separateS1S2}(f^m)$  \Comment{Reinterprets $f^m$ as $4-$dimensional}\\
\Return{${\mathfrak{J}f}$}
\end{algorithmic}
\end{algorithm}

\section{INVERSION OF 3D DRT AND 3D DJT}
\label{sec:inversion}

\subsection{Radon and John 3D adjoint transforms}
To achieve an inverse transform, it is necessary to first define an adjoint transform. If the direct transform computes the summation of values through a set of lines or planes from an input, which is complete in slopes and displacements, the adjoint transform, for every datum in the Radon or John domain corresponding to certain slopes and displacements, redistributes the energy uniformly back to the input domain by assigning a similar quantity to every voxel traversed by the line or plane. That is clearly not the inverse transform; however --and this is the key contribution from Press\cite{Press}-- it is a sufficiently good starting point to induce an iterative refinement process.

{\color{blue} These adjoint transforms will begin at stage $n$ --we will denote the backprojection partial stages as $\stackrel{\scriptscriptstyle \backsim}{f}\!{}^{m}$--, and gradually reverse the effects of the direct transform, achieving the result at $\stackrel{\scriptscriptstyle \backsim}{f}\!{}^{0}$. The mapping equations relate the same data between stages as in direct transforms, but instead of expressing where data from stage $m+1$ take their summands from, we will now express where data from stage $m$ spread their value in the direct transform; thus, the backprojection algorithm can now go to fetch them from the stage after to the previous one . Currently, the values are not the same as the original values, but are already accumulated; therefore, there is still need to perform some deconvolution or filtering to restore the original signal from $\stackrel{\scriptscriptstyle \backsim}{f}\!{}^{0}$.

The adjoint 3D DRT is defined by the inverse mapping of eq. (\ref{eq:map3D}):
\begin{equation}
\begin{split}
\label{eq:invmap3D}
\stackrel{\scriptscriptstyle \backsim}{f}\!{}^{m} ( & 
\boldsymbol{\sigma_1} \, | \, 0, \mathbf{v_1} \, | \, \boldsymbol{\sigma_2} \, | \, 0, \mathbf{v_2} \, |\,  d) 
 \,+=\, \stackrel{\scriptscriptstyle \backsim}{f}\!{}^{m+1} (  s_{1_{n-m-1}}, {\boldsymbol{\sigma_1}} |\; {\mathbf{v_1}} | \; s_{2_{n-m-1}}, {\boldsymbol{\sigma_2}} | \; {\mathbf{v_2}} | \; d)\\
\stackrel{\scriptscriptstyle \backsim}{f}\!{}^{m} (& 
\boldsymbol{\sigma_1} \, | \, 1, \mathbf{v_1} \, | \, \boldsymbol{\sigma_2} \, | \, 0, \mathbf{v_2}) \, | \, d + s_{1_{n-m-1}} + \lambda(\boldsymbol{\sigma_1})) \,+=\, \stackrel{\scriptscriptstyle \backsim}{f}\!{}^{m+1} (  s_{1_{n-m-1}}, {\boldsymbol{\sigma_1}} |\; {\mathbf{v_1}} | \; s_{2_{n-m-1}}, {\boldsymbol{\sigma_2}} | \; {\mathbf{v_2}} | \; d)\\
\stackrel{\scriptscriptstyle \backsim}{f}\!{}^{m} (
& \boldsymbol{\sigma_1} \, | \, 0, \mathbf{v_1} \, | \, \boldsymbol{\sigma_2} \, | \, 1, \mathbf{v_2} \, | \, d + s_{2_{n-m-1}} + \lambda(\boldsymbol{\sigma_2}))  \,+=\, \stackrel{\scriptscriptstyle \backsim}{f}\!{}^{m+1} (  s_{1_{n-m-1}}, {\boldsymbol{\sigma_1}} |\; {\mathbf{v_1}} | \; s_{2_{n-m-1}}, {\boldsymbol{\sigma_2}} | \; {\mathbf{v_2}} | \; d)\\
\stackrel{\scriptscriptstyle \backsim}{f}\!{}^{m} (
& \boldsymbol{\sigma_1} | 1, \mathbf{v_1} | \boldsymbol{\sigma_2} | 1, \mathbf{v_2} | 
\, d + s_{1_{n-m-1}} + \lambda(\boldsymbol{\sigma_1}) + s_{2_{n-m-1}} + \lambda(\boldsymbol{\sigma_2}))  \,+=\, \\
& \stackrel{\scriptscriptstyle \backsim}{f}\!{}^{m+1} ( s_{1_{n-m-1}}, {\boldsymbol{\sigma_1}} |\; {\mathbf{v_1}} | \; s_{2_{n-m-1}}, {\boldsymbol{\sigma_2}} | \; {\mathbf{v_2}} | \; d)
\end{split}
\end{equation}
or equivalently,
\begin{multline}
\label{eq:invmap3Dsummation}
\stackrel{\scriptscriptstyle \backsim}{f}\!{}^{m} (\boldsymbol{\sigma_1} \, | \, v_{1_m}, \mathbf{v_1} \, | \, 
\boldsymbol{\sigma_2} \, | \, v_{2_m}, \mathbf{v_2} \, |\,  d) = 
\sum_{s_{1_{n-m-1}}}^{0,\,1}\sum_{s_{2_{n-m-1}}}^{0,\,1}\stackrel{\scriptscriptstyle \backsim}{f}\!{}^{m+1}
\Big(  s_{1_{n-m-1}}, \boldsymbol{\sigma_1} |\; \mathbf{v_1} | \; 
   s_{2_{n-m-1}}, \boldsymbol{\sigma_2} | \; \\
   d \,-\,v_{1_m}\cdot\big( s_{1_{n-m-1}} + \lambda(\boldsymbol{\sigma_1})\big) \,-\, 
   v_{2_m}\cdot\big( s_{2_{n-m-1}} + \lambda(\boldsymbol{\sigma_2})\big)\Big).
\end{multline} 
Similarly, the adjoint of 3D DJT, relies on the inverse mapping of eq. (\ref{eq:map3DJ}):
\begin{multline}
\label{eq:invmap3DJsummation}
\stackrel{\scriptscriptstyle \backsim}{f}\!{}^{m} (v_m, \, \mathbf{v} \,| \, \boldsymbol{\sigma_1} \, | \, \boldsymbol{\sigma_2} \, |\, d_1 \, | \, d_2) \, = 
\sum_{s_{1_{n-m-1}}}^{0,\,1}\sum_{s_{2_{n-m-1}}}^{0,\,1}
\stackrel{\scriptscriptstyle \backsim}{f}\!{}^{m+1}
\Big({\mathbf{v}} | \; s_{1_{n-m-1}}, {\boldsymbol{\sigma_1}} | \; s_{2_{n-m-1}}, {\boldsymbol{\sigma_2}} |\; \\
d_1 - \, v_m\cdot\big(s_{1_{n-m-1}}\,+\,\lambda(\boldsymbol{\sigma_1})\big)\,|\;
d_2 - \, v_m\cdot\big(s_{2_{n-m-1}}\,+\,\lambda(\boldsymbol{\sigma_2})\big)\Big) 
\end{multline}
}


The basic backprojection algorithm for a dodecant, analogous to basic algorithm (\ref{alg:3DDRT}), is omitted because it is a reversal in the order of the $m$ loop and replacement of DRT direct formulae \ref{eq:map3D} by formulae \ref{eq:invmap3Dsummation}. 

\begin{figure}[h]
\centering
\includegraphics[width=.5\textwidth]{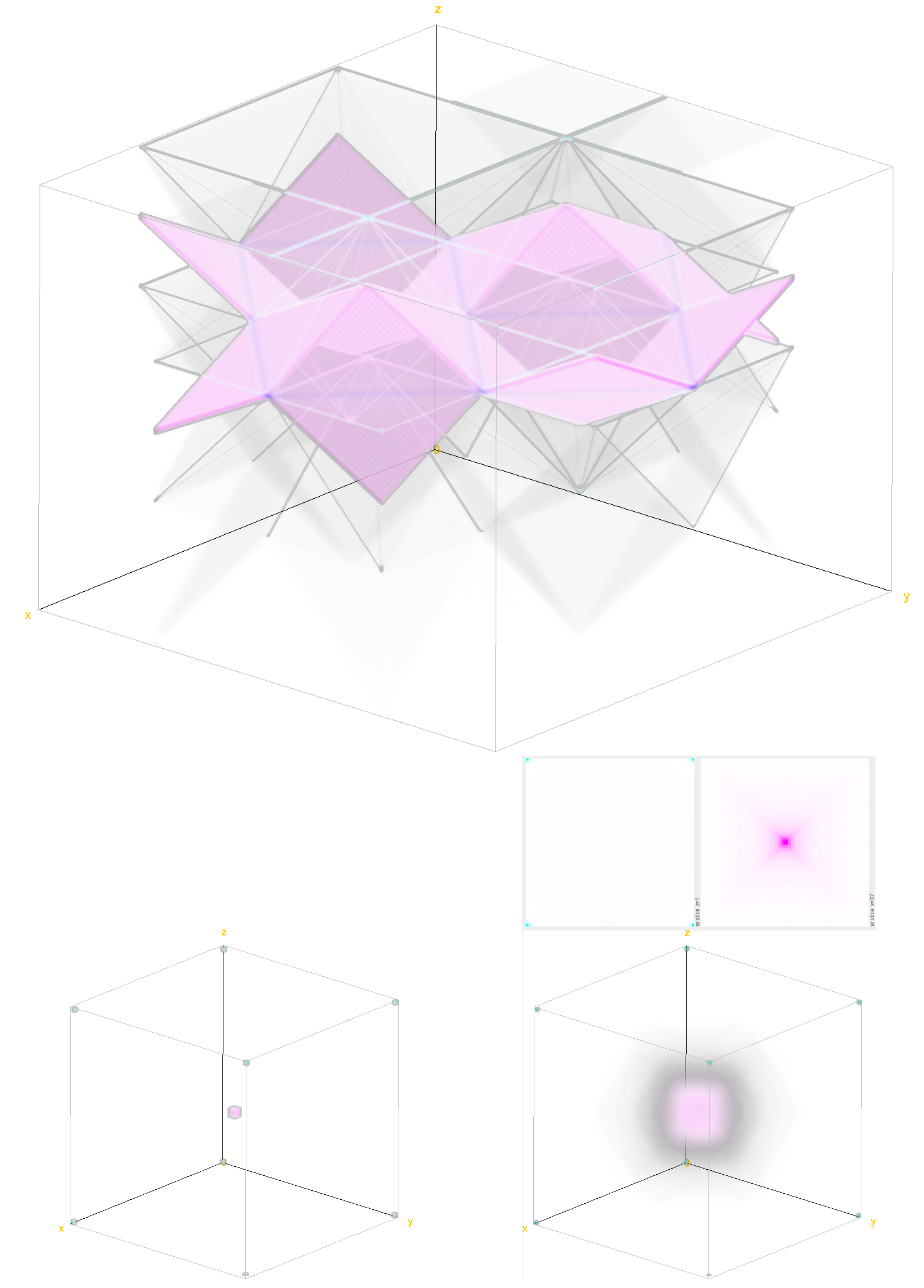}
\caption{Bottom left: Input comprising non-null values on the corner voxels of a $64^3$ cube (cyan) and the center (magenta). Top: Its 3D DRT represented as a semitransparent volume. Bottom right: ``Unfiltered'' backprojection with two slices.}
\label{fig:exampleOutput}
\end{figure}

An example of an input, its 3D DRT, and its backprojection, is presented in figure \ref{fig:exampleOutput}. An example of backprojection for the 3D DJT was already shown in figure \ref{fig:output3DDJT}. Volume visualizations have been created with \emph{ImageJ 3D viewer} \cite{imageJ3Dviewer}.

\subsection{Radon and John 3D inverse transforms}

{\color{blue} To invert these transforms we will proceed exactly as described by Press\cite{Press} for 2D, but now 3D inversion algorithms will emerge from these 3D adjoint transforms.} The adjoint transforms will be used as an approximate inverse operator of the DRT in the theory of iterative improvement of a solution to linear equations as described in chapter 2.5 of Press book on numerical methods\cite{press1982numerical}. The \emph{na\"ive} method by itself is extremely slow to converge. To accelerate this inversion so that it becomes fast, exact, and practical, Press proposes the use of a multigrid approach, described in the same book, chapter 19.6. 
The \emph{prolongation}-\emph{restriction}-\emph{filter} scheme described for two dimensions is also suitable for 3D by evolving the high-pass filter to be:
{\small
\begin{align*}
h(0,:,:)\, =\,h(2,:,:)=&
\begin{bmatrix}
-1/64 & -1/32 & -1/64\\ 
-1/32 & -1/16 & -1/32\\
-1/64 & -1/32 & -1/64
\end{bmatrix}\;\\
h(1,:,:) = &
\begin{bmatrix}
-1/32 & -1/16 & -1/32\\
-1/16 & 7/8 & -1/16\\
-1/32 & -1/16 & -1/32
\end{bmatrix}
\end{align*}}
{\color{blue} The prolongation operator $\mathcal{P}$ to be applied in the image domain:
\begin{multline}
f_{i,j,k}^{N}=f_{i+1,j,k}^{N}=f_{i+1,j+1,k}^{N}=f_{i+1,j+1,k+1}^{N}=f_{i,j+1,k}^{N}=f_{i,j+1,k+1}^{N}=f_{i,j,k+1}^{N}=f_{i+1,j,k+1}^{N}= \\ 
f_{\lfloor i / 2\rfloor, \lfloor j / 2\rfloor, \lfloor k/2 \rfloor}^{N / 2} 
\end{multline}
And the restriction operator $\mathcal{S}$ to be applied to each dodecant in the transformed domain:
\begin{equation}
\mathfrak{R}f^{N / 2}(s_1, s_2, d)=\frac{1}{8}\left[\mathfrak{R}f^{N}(2 s_1, 2 s_2, 2 d)+\mathfrak{R}f^{N}(2 s_1,2 s_2, 2 d +1)\right] 
\end{equation}
}

\begin{figure}[h]
\centering
\includegraphics[width=.6\textwidth]{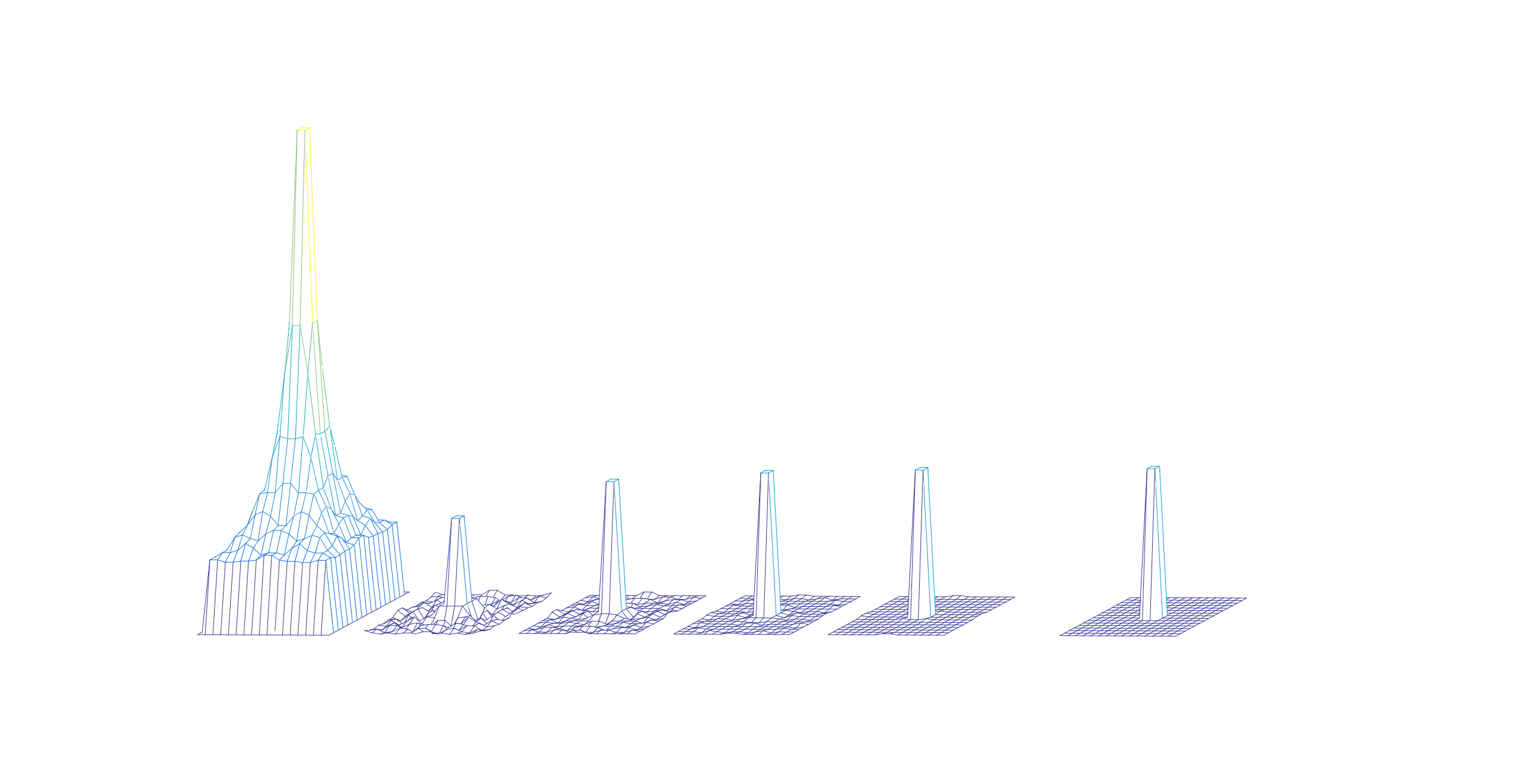}
\caption{\color{blue} Central slice profiles of inversion of input in figure \ref{fig:exampleOutput}. From left to right: backprojection; first iteration of multigrid method; fifth iteration; tenth iteration; fiftenth iteration; and input.}
\label{fig:reconstructions}
\end{figure}

The speed of convergence varies with the size $N$ and the input data; for $N=64$ about a half-dozen iterations provides very good results, and 10 iterations are almost exact (see figure \ref{fig:reconstructions}).

\begin{figure}[h]
\centering
\includegraphics[width=.65\textwidth]{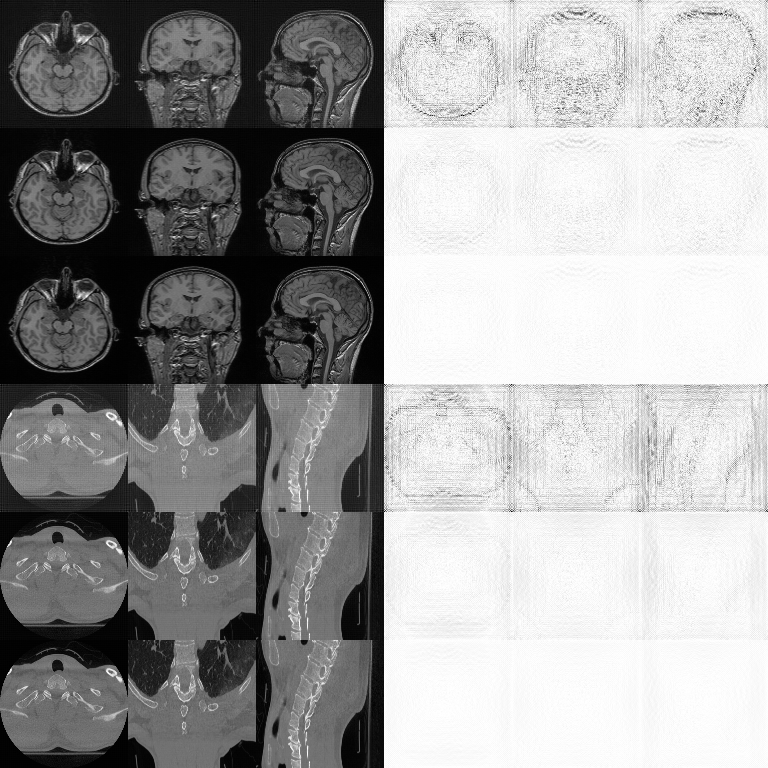}
\caption{\color{blue} Left, slices along each axis of DRT 3D inversion of a head and a spine MRI for $N = 128$. Right, absolute error. Top to bottom: inversion at iterations 5, 20 and 40.}
\label{fig:headSpine128}
\end{figure}

\begin{figure}[h]
\centering
\includegraphics[width=.5\textwidth]{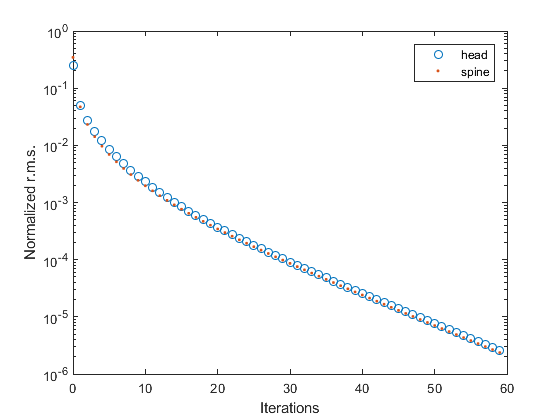}
\caption{\color{blue} Normalized root-mean-square deviation of head and spine MRI inversion.}
\label{fig:errorEvolution}
\end{figure}

{\color{blue}
Figure \ref{fig:headSpine128} shows the inversion of two MRI for $N = 128$.  The evolution of error, expressed as normalized root-mean-square deviation, is shown in figure \ref{fig:errorEvolution}, and it is very similar to that of 2D DRT, reported by Press.
In this study, we have just proved that inversion is possible by applying our forward and adjoint transforms using the method proposed by Press.} But the convergence speed of inversion could most probably be reduced with further studies on the appropriate filters and \emph{prolongation-restriction} operands. Other authors have proposed conjugate gradient methods to invert Radon transforms from the forward-backproject pair of operands.

\section{IMPLEMENTATION DETAILS AND COMPARISON}
\label{sec:implementationDetails}

\begin{figure}
\centering
\includegraphics[width=.95\textwidth]{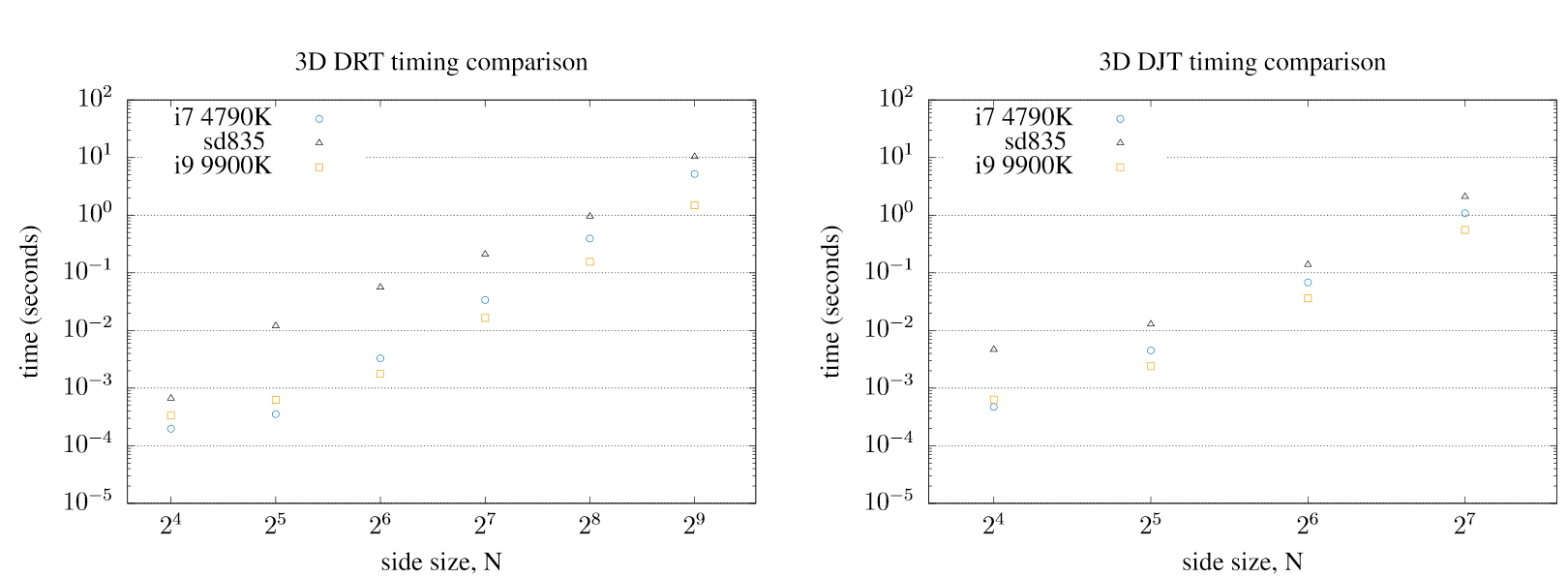}
\caption{Execution times of the 3D DRT and 3D DJT for a single dodecant on three platforms: two Intel desktop PC CPUs and a Qualcomm Snapdragon mobile phone CPU.}
\label{3D_DRT_DJT_TIMES}
\end{figure}

The 3D DRT and 3D DJT have been coded using the Halide language \cite{ragan2013halide} to obtain a reusable implementation capable of being executed in several architectures. Halide's auto-scheduling capability was used to obtain a reasonably good estimation of times without manually fine-tuning the parallel resources of each platform. Manual scheduling should reduce execution times.

Figures \ref{3D_DRT_DJT_TIMES} show the computation times required to perform a forward transform of a dodecant for different CPUs. For example, the complete set of planes traversing a dodecant of normals on a volume of $128^3$ voxels (i.e. 2 megavoxels and more than 4 million planes) can be computed in 16 milliseconds on an i9 CPU at a rate of $2.53\text{e}^8$ planes per second, where each comprises the summation of $16k$ voxels. Hardware monitoring reveals around 80\% usage of the total capacity of the i9 CPU. 

\begin{table}[]
\centering
\label{tab:condNumber}
\begin{tabular}{ccc}
 \textbf{N} & DRT3D & PPFFT3D \\
 \hline
 4 & 10.7983 & 5.8997\\
 8 & 24.1062 & 13.1113 \\ 
 16 & 67.1549 & 37.7967 
\end{tabular}
\caption{\color{blue} $L_2$ norm condition numbers for DRT and Pseudo Polar Fast Fourier based Radon transform}
\end{table}

{\color{blue} The time required to perform a backprojection is equal to that required to perform a direct transformation. However, the inversion requires the application of several backprojections (at different scales) per iteration, and subsequently several iterations; thus, the speed benefits of the proposed method in the direct path are not clear in the inverse problem. Moreover, the necessity to perform filtering destroys, as well, the ability to operate exclusively on integer arithmetics. Therefore, if the inverse transform is required, the Fourier-based methods are still a better option as the condition number  of the DRT operators is worse than that of Fourier-based operators as can be observed in Table \ref{tab:condNumber}.}

\begin{figure}[h]
\centering
\includegraphics[width=.45\textwidth]{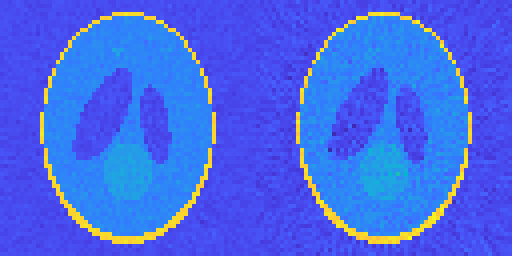}
\caption{\color{blue} A slice of the inversion of a 3D Shepp-Logan phantom of size N=64, after 60 iterations. Random noise of mean zero and variance equal to 1 percent of the coefficients magnitude has been added in the transformed domain. Left, 3D DRT inversion showing a PSNR of -15.17. Right, Averbuch's Fourier based inversion showing a PSNR of -20.02}
\label{fig:psnr-15-20}
\end{figure}
{\color{blue} However, when noise is added in the same quantity to our 3D DRT and to Averbuch's Fourier based Radon transform, the inversion is closer to the reference in ours, probably due to the lesser sparsity of the coefficientes. See figure \ref{fig:psnr-15-20}.}

{ \color{blue}
\subsection{Application example}
The speed of the proposed forward transforms operating with integers gives rise to the following application as a demonstration of its capabilities, and which distinguishes it from the rest of the existing ones in the literature, as it is capable of performing the necessary computation in video acquisition time.

\begin{figure}[h]
\centering
\includegraphics[width=.75\textwidth]{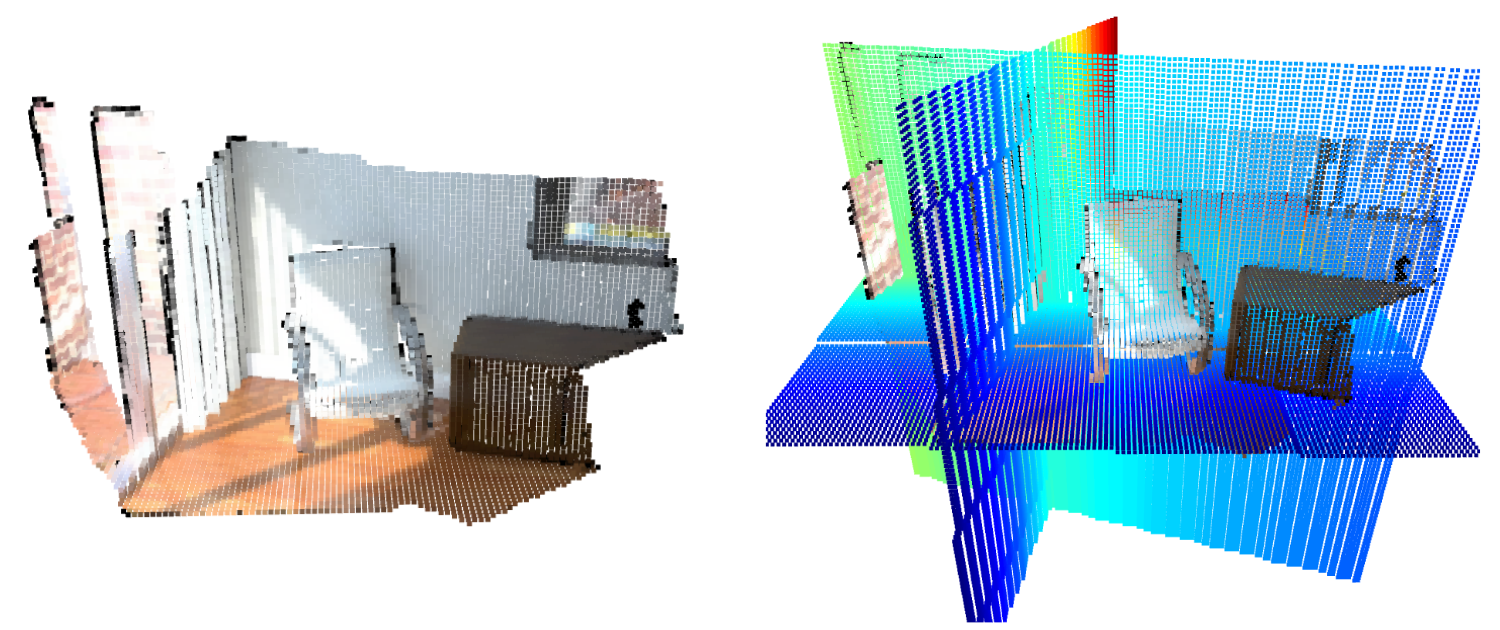}
\caption{\color{blue} Left, depth map of a scene. Right, principal planes detected by the 3D DRT of planes.}
\label{fig:applicationExample}
\end{figure}

Figure \ref{fig:applicationExample} shows an indoor scene taken with a distance sensor. This distance map is voxelized to be considered a volume cube of dimensions 128x128x128, the 3D DRT forward transform is performed and the plane that accumulated the greatest number of voxels is detected in the transformed space, and then those that are orthogonal to it, ans simultaneously are relative maxima and exceed a threshold. When those planes are backprojected, they correspond to the walls --and eventually floor and ceil--, of the room, without the need to eliminate objects or preprocess the scene. }


\subsection{Comparison with other methods}
To our knowledge, the studies conducted by Gil Shabat \emph{et al.}\cite{inverse3DFourierRadon, 3DfourierSoftware} represent the fastest implementation available for a discrete 3D Radon transform. Its main advantage is the provision of a direct inversion rather than an iterative inversion. Meanwhile, its complexity, which is asymptotically linearithmic, is similar to that of the proposed method, as well as many other Fourier-based DRTs, but its runtime complexity experiences high multiplicative factors due to interpolation. These multiplicative factors are discarded in $O()$ notation, but they make a difference when confronting execution runtimes. 

Moreover, owing to the ability of our proposed algorithms to operate in integer arithmetics, our proposed algorithms benefit from vectorized integer operations, translating into noticeable speedups in actual implementations.{\color{blue} For example, the proposed application example for detecting walls in a room could be conducted on unsigned 8 bits integers packed into 64 bytes words to obtain an immediate $64\times$ times speedup compared to other methods demanding floating point precision. In addition, we can process our 12 dodecants independently to obtain another source of immediate speedup in a 16-threads platform.}

As a result, executing on an i9 PC, Shabat's implementation of the 3D Radon transform, through the pseudo-polar Fourier domain takes 447 s to fully process a DRT forward transform on a $512^3$ input. The preprocessing times were not considered. In the same PC, it took our implementation approximately 5\% of that time: at a rate of 1.5 s per dodecant. That is, consider that even when running our implementation on a mobile phone’s CPU comprising 4-cores at 2.45 GHz, but thermally capped, it remains several times faster than Shabat's implementation running on an 8-core, 16-thread, unlocked desktop PC at 5 GHz. This advantage is lost in the inverse path.

{\color{blue} Other recent works in Radon transforms using the Fourier slice-projection theorem \cite{abbas2017exact, convolution} put their effort into eliminating or alleviating the computational cost of interpolation from Cartesian to their specially designed grids: a problem that we can simply ignore in multiscale methods.}

Another family of DRT, known as Mojettes, is based on integer arithmetics \cite{exactNTT, Mojette, fastNTFRT}. 
This approach can be extended to more than two dimensions \cite{ndimensionalMojette},
but the 'projection structures' are not lines and planes in the usual sense. Instead, they are defined with modulo arithmetics and exhibit wrap-around behavior, making them unfeasible for the line recovery example in figure \ref{fig:output3DDJT}.
Moreover, they are much slower than Fourier-based methods.

We have found two previous studies that consider the multiscale DRT extension to 3D.
Wu and Brady\cite{WuBrady} anticipated a 3D DRT of planes algorithm, not by explicitly constructing the 3D structures, but by twice applying the 2D DRT of lines, where the line integrals become plane integrals. However, they did not provide the details on the construction of the solution for a full hemisphere of vector normals.
Rim \cite{rim2017dimensional} has already tried extending multiscale DRT to 3D, and then, invert it using conjugate gradients. However, no results are provided for 3D apart from the formulae taken from Marichal \emph{et al.}\cite{Marichal}. In addition, Rim incorrectly stated that the basic algorithm would solve 1 of 16 hexadecants of normals, owing to a wrong assumption on how the four bidimensional quadrants would escalate to the third dimension.

\section{CONCLUSIONS}
\label{sec:conclusions}

Two families of algorithms, each comprising forward, adjoint, and backward 3D transforms, have been described in this paper: for Radon and John (or X-Ray) transforms. We have demonstrated the viability of the multiscale discrete approach to escalate to dimensions greater than two. These algorithms are significantly faster, in the direct path, than those already described in the literature as they elude the use of the Fourier slice-projection theorem.

Owing to their reduced computational load, they may lead the way for applications of the transforms to problems where, even if there was a solution through Radon or John transforms, the expense rendered those solutions unrealistic. In the future, we will explore the viability of the proposed algorithms in the field of computational photography \cite{radonPlenoptic, 2019Georgiev, objectTracking} and acoustic source localization.

{\color{blue} To extend the usefulness of these algorithms, we will study how to apply them to non-square dimensions. We are also exploring if a non iterative inversion is possible. }

\subsection*{Disclosures}
The authors declare that they have no conflict of interest.
\subsection*{Acknowledgments}
This work has been partially supported by 
Government of the Canary Islands and Loro Parque through "CanBio: Red de monitorización del cambio climático, la acidificación oceánica y el ruido submarino en Canarias; 2019-2022"; also by Institute for Information \& Communications Technology Promotion (IITP) grant funded by the Korean government (MSIT), [2016-0-00009], Authoring Platform Technology for Next-Generation Plenoptic Contents.
We would like to thank Editage (www.editage.com) for English language editing; and to the anonymous reviewers for their insightful comments.
\subsection*{Data, materials and code availability}
Data and code are available from the authors upon reasonable request.

\subsection*{Biographies}
\begin{wrapfigure}[5]{L}{25mm}
    \vspace{-5mm}
    \includegraphics[width=25mm,height=25mm,keepaspectratio]{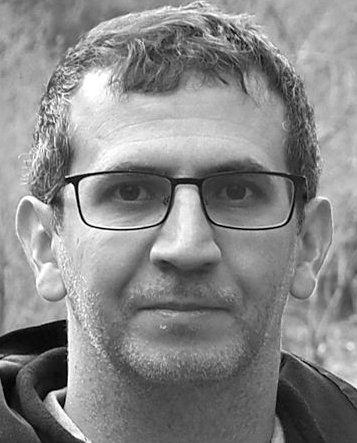}
\end{wrapfigure}
\textbf{Jos\'e G. Marichal-Hern\'andez}
  studied Computer Science in University of La Laguna, Spain. Where he also obtained his doctorate and is currently an assistant professor of Signal Theory. He is the co-inventor of three patents in effect worldwide related to 3D cameras, author of some 40 contributions to journal articles and international conferences and co-founder of the start-up Wooptix in which Intel Venture Capital has a stake.
  \par
  Other author biographies are not available.

 \bibliography{bibliografia} 
 \bibliographystyle{spiejour}

\end{document}